\documentclass[11pt]{article}
\usepackage{mathrsfs}
\usepackage{amsfonts}
\usepackage{amsfonts,latexsym,rawfonts,amsmath,amssymb,amsthm,epsfig}
\newcommand{\bthm}[2]{\vskip 8pt\noindent\bf #1\hskip 2pt\bf#2\it \hskip 8pt}
\newcommand{\ethm}{\vskip 8pt\rm}
\numberwithin{equation}{section}
\newtheorem{theorem}{Theorem}[section]
\newtheorem{lem}[theorem]{Lemma}
\newtheorem{thm}[theorem]{Theorem}
\newtheorem{pro}[theorem]{Proposition}
\newtheorem{cor}[theorem]{Corollary}

\def\s{\,\,\,\,}
\def\lan{\langle}
\def\ran{\rangle}
\def\dint{\displaystyle{\int}}
\def\dis{\displaystyle}
\def\mv{2.0ex}
\def\endproof{$\hfill\Box$\\}
\title{\bf A weak energy identity and the length of necks for a
Sacks-Uhlenbeck $\alpha$-harmonic map sequence}
\author{Yuxiang Li\thanks{This paper was written while the first author
was researching at Mathematisches Institut,
Albert-Ludwigs-Universit\"at Freiburg, supported by Alexander von
Humboldt Foundation.}\\ \and  Youde Wang\thanks{Partially supported
by 973 project of China, Grant No. 2006CB805902.} }

\date{}
\begin{document}
\maketitle

\begin{abstract}
Assume that $M$ is a closed surface and $N$ is a compact Riemannian
manifold without boundary. Let $u_\alpha: M\rightarrow N$ be the
critical point of $E_\alpha$ with $E_\alpha(u_\alpha) <C$. Assume
$u_0$ is the weak limit of $u_\alpha$ in $W^{1,2}(M, N)$ and $x_1$
is the only blow-up point in $B_\sigma(x_1)\subset M$ with $n_0$
bubbles. Then, on a local coordinate system on $B_\sigma(x_1)$ which
origin is $x_1$, we can find sequences $x_\alpha^i\rightarrow 0$,
$\lambda_\alpha^i\rightarrow 0$ ($i=1,\cdots,n_0$) s.t.
$u_\alpha(x_\alpha^i+\lambda_\alpha^ix)\rightarrow v^i$, where $v^i$
are harmonic maps from $S^2$ to $N$. We define
$$\mu_i=\liminf_{\alpha\rightarrow 1}(\lambda_\alpha^i)^{2-2\alpha}.$$
We will prove  that
$$\lim_{\alpha\rightarrow 1}E_\alpha(u_\alpha,B_\sigma(x_1))
=E(u_0,B_\sigma(x_1))+|B_\sigma(x_1)|+\sum_{j=1}^{n_0}\mu_j^2E(v^j).$$
Further,  when $n_0=1$, we define
$$\nu^1=\liminf_{\alpha\rightarrow 1}(\lambda_\alpha^1)^{-\sqrt{\alpha-1}},$$
then we have:

If $\nu^1=1$, then $u_0(B_\sigma(x_1))\cup v^1(S^2)$ is connected;

If $1<\nu^1<+\infty$, then $u_0(B_\sigma(x_1))$ and $v^1(S^2)$ are
connected by a geodesic with length
$$L = \sqrt{\dis\frac{E(v^1)}{\pi}}\log\nu^1.$$

If $\nu^1=+\infty$, the neck contains at least one geodesic
with infinite length.

\noindent We also give an example of neck which shows the neck
contains at least one geodesic of infinite length.
\\

\end{abstract}

\section{Introduction}
Let $(M,g)$ be a smooth closed Riemann surface, and $(N,h) \subset
\mathbb{R}^K$ be an n-dimensional smooth compact Riemannian
submanifold. We always assume that $N\hookrightarrow\mathbb{R}^K$ is
an isometric embedding and has no boundary.

Let $W^{1,2}(M,N)$ denote the Sobolev space of $W^{1,2}$ maps from
$M$ into $N$. If $u\in W^{1,2}(M,N)$, locally, we define the energy
density $e(u)$ of $u$ at $x\in M$ by
$$e(u)(x)=|\nabla_g u|^2=g^{ij}(x)h_{\alpha\beta}(u(x))
\frac{\partial u^\alpha}{\partial x^i}\frac{\partial
u^\beta}{\partial x^j}.$$ It is easy to check that
$$e(u)=Trace_gu^*h,$$
where $u^*h$ is the pull-back of the metric tensor $h$. Usually, the
energy $E(u)$ of $u$ is defined by
$$E(u)=\dint_Me(u)dV_g,$$
and the critical points of $E$ are called harmonic maps. We know
that a harmonic map $u$ satisfies the following equation:
$$\tau(u)=\Delta u+A(u)(\nabla u,\nabla u)=0,$$
where $A$ is the second fundamental form of $N$ in $\mathbb{R}^K$.

It is not easy to find a harmonic map, since $E$ does not satisfy
the Palais-Smale condition when the dimensions of domain manifold
$\dim(M)\ge 2$. Eells and Sampson first employed the heat flow
method to approach the existence problems of harmonic maps and
successfully deformed a map from a closed manifold into a manifold
with nonpositive sectional curvature into a homotopic harmonic map.
Concretely, they considered the heat flow for harmonic maps (or the
negative gradient flow of the energy functional $E(u)$):
$$\frac{\partial u}{\partial t}=\tau_g(u).$$
If we can establish the global existence of the above flow with
respect to the time variable $t$, or roughly speaking, the flow
flows to infinity smoothly, then we are able to find a sequence
$u_k=u(x,t_k)$ s.t. $t_k \rightarrow +\infty$ and $u_k$ converges to
a harmonic map (see \cite{E-S}).

As $\dim(M)=2$, it is well-known that the energy functional is of
conformal invariance and harmonic maps for this case are of special
importance and interest. In fact, mathematicians pay more attention
to this case. To prove the existence of harmonic maps from a closed
surface Sacks and Uhlenbeck in their pioneering paper \cite{S-U}
employed a perturbed energy functional which satisfies the
Palais-Smale condition, hence defined the so called
$\alpha$-harmonic map to approximate the harmonic map. More
precisely, for every $u\in W^{1,2\alpha}(M,N)$ Sacks and Uhlenbeck
defined the so called $\alpha$-energy $E_\alpha$ as
$$E_\alpha(u)=\dint_M(1+|\nabla u|^2)^\alpha dV_g,$$
which can be regarded as a perturbation of energy $E$, and
considered the $\alpha$-harmonic maps, i.e. the critical points of
$E_\alpha$ in $W^{1,2\alpha}(M,N)$, which satisfy the following
equation:
$$\Delta_gu_\alpha+(\alpha-1)\frac{\nabla_g|\nabla_g u_\alpha|^2
\nabla_gu_\alpha}{1+|\nabla_gu_\alpha|^2} +
A(u_\alpha)(du_\alpha,du_\alpha)=0.$$
If there is a subsequence
$u_k=u_{\alpha_k}$ which converges smoothly as $\alpha_k\rightarrow
1$, $u_{\alpha_k}$ will converge to a harmonic map.

Later, Struwe used the heat flow method of Eells and Sampson to
approach the existence problems for harmonic maps from a closed
surface and he obtained almost the same results as in \cite{S-U}.
Chang showed the same results as in \cite{St} for the case where the
domain manifold is a compact surface with smooth boundary (see
\cite{C}).

However, for both cases, the blow-up might happen. That is to say,
we are only sure that the convergence is smooth away from finitely
many points (which are called blow-up points) to a smooth harmonic
map $u_0$, which might be a trivial map. Around a blow-up point $p$,
the energy will concentrate, i.e., we will have
$$\lim_{r\rightarrow 0}\liminf_{k\rightarrow+\infty}\dint_{B_r(p)}|\nabla u_k|^2dV_g>0.$$
And then, we can find sequences
$\lim\limits_{k\rightarrow+\infty}x_k^i\rightarrow p$,
$\lim\limits_{k\rightarrow+\infty}\lambda_k^i\rightarrow 0$,
$i=1,\cdots, n_0$, s.t.
$$u_k(x_k^i+\lambda_k^ix)\rightarrow w^i\s in\s C^k_{loc}(\mathbb{R}^2 \setminus \mathcal{A}^i),$$
where all $w^i$ are non-trivial harmonic maps from $S^2$ to $N$, and
$\mathcal{A}^i$ is a finite set.

Then two problems occur. One is that if we have the energy identity, i.e.
$$\lim_{k\rightarrow+\infty}\dint_{B_\sigma}|\nabla u_k|^2dV_g=
\dint_{B_\sigma}|\nabla u_0|^2dV_g+\sum_{i=1}^{n_0}E(w^i).$$
The other
one is what the neck is if it exists?

When $u_k=u(x,t_k)$ is a subsequence of a heat flow for two
dimensional harmonic maps, the above two problems are deeply
studied. The energy identities have been proved by Qing \cite{Q} (in
the case $N=S^n$) and Ding-Tian \cite{D-T} in the general case. In
\cite{L-W}, Lin-Wang gave another proof of the energy identity. For
the neck, Qing-Tian \cite{Q-T} proved that there is no neck if the
blow-up happened at infinite time ( Ding \cite{D2} proved a more
general case), and Topping \cite{T} gave a surprising example of
heat flow blowing up at finite time s.t. the weak limit is not
continuous.

Unexpectedly, the energy identity for an $\alpha$-harmonic sequence
with bounded energy is still open. Now, many people believe that the
methods used to solve the identity for heat flow,  or more generally
a sequence with tension fields $\tau$ bounded in $L^2$, are not
powerful enough to solve the energy identity for an
$\alpha$-harmonic map sequence. The reason lies in the identity
(\ref{e1}) in this paper. For a sequence with tension fields $\tau$
bounded in $L^2$, (\ref{e1}) becomes
\begin{equation*}
    \dint_{\partial B_r}|\frac{\partial
    u_k}
      {\partial r}|^2ds_0-\frac{1}{2}\dint_{\partial B_r}
        |\nabla_0 u_k|^2ds_0=O(\dint_{B_r}|\tau(u_k)||\nabla u_k|dV_g)
    +O(1).
\end{equation*}
then the right side of the above identity is bounded. However, in
(\ref{e1}), a very bad term
$$\frac{\alpha-1}{r}\int_{B_r}(1+|\nabla u_\alpha|^2)^{\alpha-1}
|\nabla u_\alpha|^2dV_g$$
appears.

The known energy identities for some special $\alpha$-harmonic
sequences are usually obtained by methods which are completely
different with the one of \cite{D-T}. Now we would like to mention
the following cases.

If $\{u_\alpha\}$ is a sequence of minimizing $\alpha$-harmonic map,
i.e. every $u_\alpha$ is the minimizer of $E_\alpha$, which belongs
to the same homotopic class, Chen and Tian \cite{C-T} proved that
the necks consist of some geodesics of finite length, and moreover
this implies no loss of energy in necks for the sequence (see also
\cite{D-K}).

Another important case is the energy identity for a minimax
sequence.  We let $M$ be a compact Riemann surface, $A$ be a
parameter manifold. Let $h_0:M\times A\rightarrow N$ be continuous.
Assume $H$ be the class of all maps homotopic to $h_0$, and
\begin{equation}\label{jost}
\beta_\alpha(H)=\inf_{h\in H}\sup_{t\in A}E_\alpha(h(\cdot,t)).
\end{equation}
We can deduce from Jost's  result \cite{J} that  there is at least
one sequence  $u_{\alpha_k}$ which attained $\beta_{\alpha_k}(H)$
satisfies the energy identity as $\alpha_k\rightarrow 1$ (Also see
\cite{C-M} and \cite{L}).

In this paper, we will adopt some methods and techniques in
\cite{D-T} and \cite{D2} to discuss the energy identity for an
$\alpha$-harmonic sequence, especially the necks between the
bubbles. However, we can not give a final proof on the energy
identity for such Sacks-Uhlenbeck sequence, instead, we only show a
weaker energy identity and give some observation on this subject.
On the other hand, we exploit the details of the necks.
Precisely we provide a new method to show that the necks converge to geodesics
and obtain the formula on the length of the geodesics.

Now, we assume that $u_\alpha$ is a sequence of $\alpha$-harmonic
maps from $(M,g)$ to $(N,h)$ with
$$E_{\alpha}(u_\alpha)<\Theta.$$
Then, by the theory of Sacks and Uhlenbeck, we are able to assume
that there exists a sequence $\alpha_k\rightarrow 1$, s.t.
$u_{\alpha_k}$ converges to a harmonic map $u_0:M\rightarrow N$
smoothly away from a finite many points $\{x_i\}$ as
$\alpha_k\rightarrow 1$. We assume that there are $n_0$ bubbles at
the point $x_1$. Then we are able to assume that there are
$x_{\alpha_k}^j\rightarrow x_1$ and $\lambda_{\alpha_k}^j\rightarrow
0$ for $j=1,\cdots,n_0$, such that
$$v_{\alpha_k}^j=u_{\alpha_k}(x_{\alpha_k}^j+\lambda_{\alpha_k}^jx)$$
converge in
$C^k_{loc}({\mathbb{R}^2\setminus\{p_1,p_2,\cdots,p_{s_j}\}})$ to
non-trivial harmonic maps
$$v^j: S^2\rightarrow N.$$
Moreover, we assume that one of the following holds:\\

{\bf H1}. For any fixed $R$, $B_{R\lambda_{\alpha_k}^i}(x_{\alpha_k}^i)\cap B_{R\lambda_{\alpha_k}^j}
(x_{\alpha_k}^j)=\emptyset$ whenever
$({\alpha_k}-1)$ are sufficiently small.

{\bf H2}. $\frac{\lambda_{\alpha_k}^i}{\lambda_{\alpha_k}^j}+\frac{\lambda_{\alpha_k}^j}{\lambda_{\alpha_k}^i}
\rightarrow+\infty$
as ${\alpha_k}\rightarrow 1$.\\

\bthm{Remark}{1.} One is easy to check that if
$(\lambda_{\alpha_k}^i,x_{\alpha_k}^i)$ and $(\lambda_{\alpha_k}^j,x_{\alpha_k}^j)$
do not satisfy {\bf H1} and {\bf H2}, then we can find  subsequences
of $\lambda_{\alpha_k}^i$, $x_{\alpha_k}^i$ and $\lambda_{\alpha_k}^j$,
$x_{\alpha_k}^j$ s.t.
$\frac{\lambda_{\alpha_k}^i}{\lambda_{\alpha_k}^j}\rightarrow \lambda\in
(0,\infty)$ and
$\frac{x_{\alpha_k}^i-x_{\alpha_k}^j}{\lambda_{\alpha_k}^j}\rightarrow a\in
\mathbb{R}^2$. Since
$$u_{\alpha_k}(x_{\alpha_k}^i+\lambda_{\alpha_k}^i x)=u_{\alpha_k}(x_{\alpha_k}^j
+\lambda_{\alpha_k}^j(\frac{x_{\alpha_k}^i-x_{\alpha_k}^j}
{\lambda_{\alpha_k}^j}+\frac{\lambda_{\alpha_k}^i}{\lambda_{\alpha_k}^j}x)),$$
we have
$$v^i(x)=v^j(a+\lambda x),$$
and then $v^i$ and $v^j$ are in fact the same bubble. \ethm

Fixing an $R$, we have
\begin{equation}
\begin{array}{ll}
&\dint_{B_{R\lambda_{\alpha_k}^j}(x_{\alpha_k}^j)\setminus
(\cup_{i=1}^{s_j}B_{\delta\lambda_{\alpha_k}^j} (x_{\alpha_{k}}^j +
\lambda_{\alpha_k}^jp_{i}))}
|\nabla_g u_{\alpha_k}|^{2\alpha_k}dV_g \\
=&(\lambda_{\alpha_k}^j)^{2-2\alpha} \dint_{B_R\setminus
(\cup_{i=1}^{s_j}B_{\delta}(p_{i}))}|\nabla_g
v_{\alpha_k}^j|^{2\alpha_k}dV_{g(x_{\alpha_k}^j+\lambda_{\alpha_k}^jx)}.
\end{array}
\end{equation}
Since
$$\dint_{B_R\setminus (\cup_{i=1}^{s_j}B_{\delta}(p_i))}|\nabla_g v_{\alpha_k}^j|^{2\alpha_k}
dV_{g(x_{\alpha_k}^j+\lambda_{\alpha_k}^j x)} \rightarrow
\dint_{B_R\setminus (\cup_{i=1}^{s_j}B_{\delta}(p_i))}|\nabla_0
v^j|^2dx,$$ and
$$\lambda_{\alpha_k}^j<1,$$
we define
\begin{equation}\label{dm}
\begin{array}{lll}
\mu_j&=&\liminf\limits_{\alpha\rightarrow 1}
(\lambda_{\alpha}^j)^{2-2\alpha}\leq\lim\limits_{k\rightarrow
\infty}\frac{\dint_{B_{R\lambda_{\alpha_k}^j}(x_{\alpha_k}^j)\setminus
(\cup_{i=1}^{s_j}B_{\delta\lambda_{\alpha_k}^j} (x_{\alpha_k}^j +
\lambda_{\alpha_k}^j p_i))}
|\nabla_gu_{\alpha_k}^j|^{2\alpha_k}dV_g}{\dint_{B_R\setminus
(\cup_{i=1}^{s_j}B_{\delta}(p_i))}|\nabla_0 v^j|^2dx}\\
&\leq& \frac{\Theta-|M|-E(u_0)}{\theta},
\end{array}
\end{equation}
where
$$\theta=\inf\{E(u): u\hbox{ is a nontrivial harmonic map from } S^2 \hbox{ to } N\}.$$
Therefore, we know
$$\mu_j \in [1,\frac{\Theta-|M|-E(u_0)}{\theta}].$$

The first task of  this paper is to get the following weak energy
identity:

\begin{thm}\label{main}
Let $M$ be a smooth closed Riemann surface and $N$ be a smooth
compact Riemannian manifold without boundary. Assume that
$u_{\alpha_k}\in C^\infty(M, N)$ $(\alpha_k\rightarrow 1)$ is a
sequence of $\alpha_k$-harmonic maps with uniformly bounded energy
and $x_1$ be the only blow-up point of the sequence
$\{u_{\alpha_k}\}$ in $B_\sigma(x_1)\subset M$. Then, passing to a
subsequence, there exist $u_0: M\rightarrow N$ which is a smooth
harmonic map and finitely many bubbles $v_j: S^2\rightarrow N$ such
that $u_{\alpha_k}\rightarrow u_0$ weakly in $W^{1,2}(M, N)$ and in
$C^\infty_{loc}(B_\sigma(x_1) \setminus\{x_1\}, N)$ and the
following identity holds
\begin{equation}\label{ee}
\lim_{k\rightarrow +\infty}E_{\alpha_k}(u_{\alpha_k},B_\sigma(x_1))
=E(u_0,B_\sigma(x_1))+|B_\sigma(x_1)|+\sum_{j=1}^{n_0}\mu_j^2E(v^j),
\end{equation}
where $\mu_j$ is defined by (\ref{dm}) and $n_0$ is the number of
bubbles at $x_1$.
\end{thm}

This theorem tells us that the energy identity holds true if and
only if $\mu_j=1$. It provides a new route to approach the problem
whether the necks contain energy or not.

\bthm{Remark}{2.} By Lemma \ref{a1} in section 2, $\mu_j=1$ implies
\begin{equation}\label{ee1}
\lim_{k\rightarrow +\infty}E(u_{\alpha_k},B_\sigma(x_1))
=E(u_0,B_\sigma(x_1))+\sum_{j=1}^{n_0}E(v^j),
\end{equation}
and reversely, by Lemma \ref{a1} and (\ref{e11}), (\ref{ee1}) also
implies $\mu_j=1$. \ethm

It is our another purpose to study the behavior of the necks
connecting bubbles. For this sake, we need to define
$$\nu_j=\liminf_{\alpha\rightarrow 1}(\lambda_\alpha^j)^{-\sqrt{\alpha-1}}.$$
We will see that the above quantity play an important role in the
discussion on the behavior of blowing up. Our main results are
stated as follows:

\begin{thm}\label{main2}
Let $M$ be a smooth closed Riemann surface and $N$ be a smooth
closed Riemannian manifold and $u_{\alpha_k}\in C^\infty(M, N)$ be a
sequence of $\alpha_k$-harmonic maps with uniformly bounded energy
and $u_{\alpha_k}$ converges to a smooth harmonic map $u_0:
M\rightarrow N$ in $C^\infty_{loc}(B_\sigma(x_1) \setminus \{x_1\},
N)$ as $\alpha_k\rightarrow1$. Assume there is only one bubble in
$B_\sigma(x_1) \subset M$ for $\{u_{\alpha_k}\}$ and $v^1:
S^2\rightarrow N$ is the bubbling map. Let
$\nu^1=\liminf\limits_{\alpha\rightarrow
1}(\lambda_\alpha^1)^{-\sqrt{\alpha-1}}$. Then we have

1) when $\nu^1=1$, the set $u_0(B_\sigma(x_1))\cup v^1(S^2)$ is a
connected subset of $N$;

2) when $\nu^1\in(1,\infty)$, the set $u_0(B_\sigma(x_1))$ and
$v^1(S^2)$ are connected by a geodesic with Length $$L =
\sqrt{\frac{E(v^1)}{\pi}}\log\nu^1;$$

3) when $\nu^1=+\infty$, the neck contains at least an infinite length geodesic.
\end{thm}

\bthm{Remark}{3.} Although we state and prove Theorem \ref{main2}
only for one bubble case, it is not difficult to follow the steps in
section 3.2 to prove the general case. However,  the general case is
quite complicated, for example, if we have 2 bubbles:
$$u_\alpha(\lambda_\alpha^1x+x_1)\rightarrow v^1,\s and \s
u_\alpha(\lambda_\alpha^2x+x_1)\rightarrow v^2$$
which satisfy: $\lambda_\alpha^1/\lambda_\alpha^2\rightarrow 0$ and
$\nu^1,\nu^2<\infty$, then
 $u_0(B_\delta(x_1))$,
$v^2(S^2)$ are connected by a geodesic with length
$$L=\sqrt{\frac{E(v^1)+E(v^2)}{\pi}}\log\nu^2,$$
and
 $v^1(S^2)$, $v^2(S^2)$ are connected by a geodesic with length
$$L=\sqrt{\frac{E(v^1)}{\pi}}\log\frac{\nu^1}{\nu^2}.$$\ethm

We should mention that after we completed the paper we found that
Moore had proved that if a neck is of finite length $L$ and
$\tilde{g}\geq 1$ (the genus of $M$), then $L =
\sqrt{\frac{E(v^1)}{\pi}} \log\nu$ (note that in \cite{M}, $E(u)$ is
defined to be $\frac{1}{2}\int_M|\nabla u|^2dV_g$). However, the
arguments to prove Theorem \ref{main2} in this paper is completely
different from Moore's proof. The key estimation of us is the
Proposition \ref{p1} in section 4, which gives the details of the
necks.

The Proposition \ref{p1} also provides a new method to prove that
the necks consist of geodesics, which has been already proved by
Chen and Tian \cite{C-T}. In this paper, we will make use of the
following curve
$$\Gamma_\alpha(r)=\frac{1}{2\pi}\int_0^{2\pi}u_\alpha(r,\theta)d\theta$$
to approximate the necks. With the help of Proposition \ref{p1}, one
can easily calculate the second fundamental form of the
approximation curve, and then to prove that the limiting curve
satisfies the equation of geodesic in $N$.

We failed  to find a  sufficient condition s.t. $\nu^i<+\infty$, but
we will show that there are indeed many cases that the necks contain
at least one infinite length geodesic:

\begin{cor}\label{intro.c1} Let $\alpha_k\rightarrow 1$, and $u_{k}:
M\rightarrow N$ be a minimizer of $E_{\alpha_k}$ in the homotopic
class containing $u_k$. We assume for any $i\neq j$, $u_{i}$ and
$u_{j}$ are not in the same homotopic class. If
$$\sup_{k}E_{\alpha_k}(u_k)<+\infty,$$
then $u_k$ will blow up, and the neck contains at least one infinite
length geodesic.
\end{cor}

\bthm{Remark}{3.}In the last section, by constructing a manifold $N$
we will give an example of a minimizing $\alpha$-harmonic map
sequence, which satisfies the condition in the above corollary. This
indicates that there exists a neck joining bubbles which is a
geodesic of infinite length. \ethm

\noindent We conclude this introduction with showing the following
proposition as a consequence of Theorem \ref{main}, which implies the
result due to Chen-Tian that, if the necks consist of some geodesics
of finite length, then the energy identity is true:

\begin{pro} The energy identity holds true for a subsequence of $u_{\alpha}$ if and only if
\begin{equation}\label{in1}
\liminf_{\alpha\rightarrow1}\|\nabla u_{\alpha}\|_{C^0(M)}^{\alpha-1}=1.
\end{equation}
The limit set of such subsequence has no neck if and only if
\begin{equation*}
\liminf_{\alpha\rightarrow1}\|\nabla u_{\alpha}\|_{C^0(M)}^{\sqrt{\alpha-1}}=1.
\end{equation*}
The bubbles in limit set of such subsequence are joined by some
geodesics of finite length, if and only if
\begin{equation*}
\liminf_{\alpha\rightarrow1}\|\nabla u_{\alpha}\|_{C^0(M)}^{\sqrt{\alpha-1}}<+\infty.
\end{equation*}
\end{pro}

\proof We only prove the first claim.

First, we prove (\ref{in1}) implies $\mu_j=1$. We assume
$v_\alpha^j(x)=u_\alpha(x_\alpha^j+\lambda_\alpha^jx)$ converges to
$v^j$ in $C^1_{loc}(\mathbb{R}^n\setminus \{p_1,p_2,\cdots,p_s\})$.
Then we have
$$(\lambda_\alpha^j)^{1-\alpha}=\frac{|\nabla u_\alpha(x_\alpha^j+\lambda_\alpha^j x)|^{\alpha-1}}
{|\nabla v_\alpha^j(x)|^{\alpha-1}}$$ for any $x$ with $|\nabla
v^j(x)|\neq 0$. Hence we get $\mu_j\leq 1$ and then $\mu_j=1$.

Now, we will prove ``$\mu_j=1$ for all $j$" implies (\ref{in1}). Let
$x_\alpha$ to be the point s.t. $|\nabla
u_\alpha|(x_\alpha)=\max|\nabla u_\alpha|$, and
$$\lambda_\alpha=\frac{1} {|\nabla u_\alpha|(x_\alpha)}.$$ We set
$v_\alpha(x)=u_\alpha(x_\alpha+\lambda_\alpha x)$. One is easy to
check that $v_\alpha$ will converge to a non-trivial harmonic map
$v_0$ locally. By {\bf H1} and {\bf H2} we must find a $j$, s.t.
$$B_{R\lambda_\alpha^j}(x_\alpha^j)\cap B_{R\lambda_\alpha}(x_\alpha)\neq\emptyset, \s and \s
\frac{1}{C}\lambda_\alpha^j<\lambda_\alpha<C\lambda_\alpha^j$$
for some $C>0$. Hence we get $|\lambda_\alpha|^{\alpha-1}\rightarrow 1$.
\endproof

{\bf Acknowledgement:} {\small The authors is grateful to thank
Professor W. Ding for his help and encouragement. The first author
would like to thank Prof. E. Kuwert for many helpful discussions.}

\section{Preliminary}
In this section we intend to establish some integral formulas on
$\alpha$-harmonic maps from a closed surfaces by the variations of
domain. Of course, we need to choose some suitable variational
vector fields on $M$ which generate the transformations of $M$. We
will see that these integral relations will play an important role
in the proofs of main theorems.

Note that the functional $E_\alpha$ is not conformal invariant. For
example, on an isothermal coordinate system around a point $p\in M$,
if we set the metric
$$g=e^{\varphi}((dx)^2+(dy)^2)$$
with $p=(0,0)$, $\varphi(0)=0$ and $\tilde{u}_\alpha(x)=u_\alpha(\lambda x)$,
then we will get
$$\int_{B_\delta}(1+|\nabla_g u_\alpha|^2)^{\alpha}dV_g=
\int_{B_{\frac{\delta}{\lambda}}}\lambda^{2-2\alpha}(\lambda^2+
|\nabla_{g'}\tilde{u}_\alpha|^2)^\alpha dV_{g'},$$
where $g'=e^{\varphi(p+\lambda x)}((dx)^2+(dy)^2)$.
We also ought to note that an $\alpha$-harmonic map sequence $u_\alpha$ may
have several bubbles near a blowing up point, for example, there are sequences
$\lambda_\alpha^1$, $\lambda_\alpha^2$, s.t.
$$\frac{\lambda_\alpha^1}{\lambda_\alpha^2}\rightarrow 0,\s \lambda_\alpha^2\rightarrow 0,$$
as $\alpha\rightarrow 1$, and
$$v_\alpha^1(x)=u_\alpha(\lambda_\alpha^1x)\rightarrow v^1\s in\s C^k_{loc}
(\mathbb{R}^2),\s
v_\alpha^2(x)=u_\alpha(\lambda_\alpha^2x)\rightarrow v^2\s in\s
C^k_{loc}(\mathbb{R}^2\setminus\{0\}),$$
where $v^1$, $v^2$ are
non-trivial harmonic maps from $S^2$ to N. For this case, we have
$$v^1_\alpha(x)=v^2_\alpha\left(\frac{\lambda_\alpha^1}{\lambda_\alpha^2}x\right),$$
i.e. $v^1(x)$ is in fact a bubble for the sequence $v_\alpha^2$.
Therefore, we need to consider the equation of $v^2_\alpha$, and one
is easy to check that $v_\alpha^2$ is locally a critical point of
the functional
$$F(v)=\int_{B_\delta}((\lambda_\alpha^2)^2+|\nabla v|^2)^\alpha dV_{g_\alpha},$$
where $g_\alpha=e^{\varphi(\lambda^2_\alpha x)}((dx)^2+(dy)^2)$. For
this reason, we need to consider a more general $\alpha$-energy
which is of the following form:
$$E_{\alpha,\epsilon_\alpha}(u)=\dint_{B_\delta}(\epsilon_\alpha+|\nabla_gu_\alpha|^2)^\alpha
 dV_{g_\alpha}.$$
Let $u_\alpha$ be the critical point of the above functional. Then,
$u_\alpha$ satisfies the following elliptic system which is also
called the equation of $\alpha$-harmonic maps:
\begin{equation}\label{eq1}
\Delta_{g_\alpha}u_\alpha+(\alpha-1)\frac{\nabla_{g_\alpha}|\nabla_{g_\alpha}u_\alpha|^2\nabla_{g_\alpha}u_\alpha}
{\epsilon_\alpha+
|\nabla_{g_\alpha}u_\alpha|^2}+A(u_\alpha)(du_\alpha,du_\alpha)=0.
\end{equation}
Here we always assume that the sequence $\epsilon_\alpha$
($\epsilon_\alpha\leq 1$) satisfies
\begin{equation}\label{beta0}
\lim\limits_{\alpha
\rightarrow 1}{\epsilon_\alpha}^{\alpha-1}> \beta_0>0.
\end{equation}
It follows
from (\ref{dm}) that this assumption is reasonable.

From now on, we consider $u_\alpha$ to a map sequence from
$(B,g)$ to $(N,h)$ which satisfy equation (\ref{eq1}). We assume that $g=e^{\varphi_\alpha}(
(dx^1)^2+(dx^2)^2)$ with $\varphi_\alpha(0)=0$
and $\varphi_\alpha\rightarrow \varphi$ smoothly.
Moreover, we assume that $u_\alpha\rightarrow u_0$
in $C^k_{loc}(\bar{B}\setminus\{0\})$.\\

\noindent Next, we recall the well-known $\epsilon$-regularity
theorem due to Sacks-Uhlenbeck \cite{S-U}:

\begin{thm}\label{epsilon}
Let $u: B\rightarrow N$ satisfies equation (\ref{eq1}) where
$B\subset M$ is a ball with radius $1$. There exists $\epsilon>0$
and $\alpha_0>1$ such that if $E(u,B)<\epsilon$ and
$1\leq\alpha\leq\alpha_0$,  then for all smaller $r<1$, we have
$$\|\nabla u\|_{W^{1,p}(B_r)}\leq C(p,r)E(u,B),$$
here $B_r\subset B$ is a ball with radius $r$, $1 < p < \infty$.
\end{thm}

\noindent We also have
\begin{lem}\label{a1} Let $u_\alpha$ be the critical point of $E_\alpha$ with $E_\alpha(u_\alpha)\leq
\Theta$. We have
$$\beta_0<\liminf_{\alpha\rightarrow 1}\|(\epsilon_\alpha+|\nabla_gu_\alpha|^2)^{\alpha-1}\|_{C^0(B)}
\leq\limsup_{\alpha\rightarrow 1}\|(\epsilon_\alpha+|\nabla_gu_\alpha|^2)^{\alpha-1}\|_{C^0(B)}<\beta_1,$$
where $\beta_1$ is independent  of $\alpha$.
\end{lem}
\proof
Obviously, we only need  to prove
$\|\nabla_gu_\alpha\|^{\alpha-1}_{C^0(B)}<C$. We assume that there is sequence $\alpha_k\rightarrow 1$, s.t.
 $\|\nabla_gu_{\alpha_k}\|^{\alpha_k-1}_{C^0(B)}
\rightarrow+\infty$ as $k\rightarrow+\infty$.

Let
$|\nabla_gu_{\alpha_k}|(x_{\alpha_k})=\max\{|\nabla_gu_{\alpha_k}|\}$,
and $\lambda_k=\frac{1}{|\nabla_gu_{\alpha_k}|}$, $v_k(x)=
u_{\alpha_k}(x_{\alpha_k}+\lambda_k x)$. Then a subsequence of
$\{v_{k_j}\}$ converges to a new nontrivial harmonic map from $S^2$
to $N$. Then by (\ref{dm}), we obtain the following
$\lambda_{k_j}^{1-\alpha_{k_j}}<C$, which contradicts with the
choice of $\alpha_k$.
\endproof

\subsection{Variational formula}
Take an $1$-parameter family of transformations $\{\phi_s\}$
which is generated by the vector field $X$. If we assume $X$ is supported in $B$, then we have
$$\begin{array}{ll}
   E_{\alpha,\epsilon_\alpha}(u\circ\phi_s)&=\dint_B(\epsilon_\alpha+|\nabla_g (u\circ\phi_s)|^2)^\alpha
   dV_g\\[\mv]
      &=\dint_B(\epsilon_\alpha+\sum_\beta|d(u\circ\phi_s)(e_\beta(x))|^2)^\alpha
      dV_g(x)\\[\mv]
      &=\dint_B(\epsilon_\alpha+\sum_\beta|du({\phi_{s}}_*(e_\beta(x)))
       |^2)^\alpha dV_g(x)\\[\mv]
     &=\dint_B(\epsilon_\alpha+\sum_\beta|du({\phi_{s}}_*(e_\beta(\phi_s^{-1}(x))))
     |^2)^\alpha
     Jac(\phi_{s}^{-1})dV_g,
  \end{array}$$
where $\{e_\alpha\}$ is a local orthonormal basis of $TB$. Noting
 $$\frac{d}{ds}Jac(\phi_s^{-1})dV_g|_{s=0}=-div(X)dV_g,$$
we have proved the formula
$$\begin{array}{ll}dE_f(u)(u_*(X)) =& -\dint_B(\epsilon_\alpha+|\nabla_g
u|^2)^\alpha div(X)dV_g\\
&+2\alpha\sum\limits_\beta\dint_B(\epsilon_\alpha+|\nabla_g
u|^2)^{\alpha-1} \lan du(\nabla_{e_\beta} X),du(e_\beta)\ran
dV_g.\end{array}$$

Now, we assume $u_\alpha$ to be the critical point of $E_\alpha$.
For any vector field $X$ on $B$, we have
$$-\dint_B(\epsilon_\alpha+|\nabla_g u_\alpha|^2)^\alpha div X dV_g+2\alpha
\sum_{\beta}\dint_B(\epsilon_\alpha+|\nabla_g
u_\alpha|^2)^{\alpha-1} \lan
du_\alpha(\nabla_{e_\beta}X),du_\alpha(e_\beta)\ran dV_g=0.$$

Next, for $0<t'<t\leq\rho$, we choose a vector field $X$ with
compact support in $B_\rho$ by $X=\eta(r)r\frac{\partial}{\partial
r}=\eta(|x|)x^i\frac{\partial}{\partial x^i}$, where $\eta$ is
defined by
$$\eta(r)=\left\{\begin{array}{ll}
                     1&if\s r\leq t'\\[\mv]
                     \dis\frac{t-r}{t-t'}&if\s t'\leq r\leq t\\[\mv]
                     0&if\s r\geq t,
                  \end{array}\right.$$
where $r=\sqrt{(x^1)^2+(x^2)^2}$. By a direct computation we obtain
$$div(X)=2\eta+r\eta'+r\eta\frac{\partial{\varphi}}{\partial r},$$
and
$$\nabla_{\frac{\partial}{\partial x^1}}X=\eta\frac{\partial}{\partial x^1}
+\eta'\frac{(x^1)^2}{r}\frac{\partial}{\partial
x^1}+\eta'\frac{x^1x^2}{r}\frac{\partial}{\partial x^2} +\eta
x^1\Gamma_{11}^1\frac{\partial}{\partial x^1} +\eta
x^1\Gamma_{11}^2\frac{\partial}{\partial x^2} \\
+ \eta x^2\Gamma_{12}^1\frac{\partial}{\partial x^1}+\eta
x^2\Gamma_{12}^2\frac{\partial} {\partial x^2}.$$
Then,
$$\begin{array}{lll}
   \sum\limits{_\beta}\lan du_\alpha(\nabla_{e_\beta}X),du_\alpha(e_\beta)\ran dV_g&=&
       \lan du_\alpha(\nabla_\frac{\partial}{\partial x^1}X),
        du_\alpha(\frac{\partial}{\partial x^1})\ran dx
       +\lan du_\alpha(\nabla_\frac{\partial}{\partial x^2}X),
       du_\alpha(\frac{\partial}{\partial x^2})
       \ran dx\\[\mv]
   &=&(\eta|\nabla_0u_\alpha|^2+\eta'r|\frac{\partial u_\alpha}{\partial r}|^2+O(|x|)
    |\nabla_0 u_\alpha|^2)dx,
  \end{array}$$
where $\nabla_0$ is the Riemannian connection with respect to
standard metric. Hence, we derive
$$\begin{array}{lll}
   0&=&(2\alpha-2)\dint_{B_t}\eta(\epsilon_\alpha+|\nabla_gu_\alpha|^2)^{\alpha-1}|\nabla_0
   u_\alpha|^2dx\\[\mv]
    &&+\dint_{B_t}O(|x|)(\epsilon_\alpha+|\nabla_gu_\alpha|^2)^{\alpha-1}|\nabla_0
    u_\alpha|^2dx\\[\mv]
     &&-2\epsilon_\alpha\dint_{B_t}\eta(\epsilon_\alpha+|\nabla_g u_\alpha|^2)^{\alpha-1}dV_g
     +\dis\frac{\epsilon_\alpha}{t-t'}\dint_{B_t\setminus B_{t'}}r(\epsilon_\alpha+|\nabla_gu|^2)^{\alpha-1}
     dV_g\\[\mv]
    &&+\dis\frac{1}{t-t'}\dint_{B_t\setminus B_{t'}}(\epsilon_\alpha+|\nabla_gu_\alpha|^2)^{\alpha-1}
     [|\nabla_0u_\alpha|^2r-2\alpha r
      |\frac{\partial u_\alpha}
       {\partial r}|^2]dx\\[\mv]
    &&-\dint_{B_t}\epsilon_\alpha(\epsilon_\alpha+|\nabla_gu_\alpha|^2)^{\alpha-1}r\eta\frac{\partial{\varphi}}
    {\partial r} dV_g.
  \end{array}$$
Letting $t'\rightarrow t$ in the above identity and using Lemma
\ref{a1}, we obtain the following
\begin{equation}\label{e1}
\begin{array}{l}
    \dint_{\partial B_t}(\epsilon_\alpha+|\nabla_g u_\alpha|^2)^{\alpha-1}|\frac{\partial
    u_\alpha}{\partial r}|^2ds_0-\frac{1}{2\alpha}\dint_{\partial B_t}(\epsilon_\alpha
    +|\nabla_g u_\alpha|^2)^{\alpha-1}|\nabla_0 u_\alpha|^2ds_0\\[\mv]
    \s\s\s\s=\dis\frac{(\alpha-1)}{\alpha t}\dint_{B_t}(\epsilon_\alpha+|\nabla_gu_\alpha|^2
    )^{\alpha-1}|\nabla_0u_\alpha|^2dx
    +O(t),
\end{array}
\end{equation}
where $ds_0$ is the volume element of $\partial B_t$ with respect to
the Euclidean metric. We know that the metric $g$ can be written as
$g=e^\varphi(dr^2 + r^2d\theta^2)$ in the polar coordinate system.
Set
 $$u_{\alpha,\theta}=\frac{1}{r}\frac{\partial
 u_\alpha}{\partial\theta}.$$
Since $|\nabla_0u_\alpha|^2=|\frac{\partial u_\alpha}{\partial
r}|^2+|u_{\alpha,\theta}|^2$, we get from the above identity
\begin{equation}\label{e11}
\begin{array}{l}
    (1-\dis\frac{1}{2\alpha})\dint_{\partial B_t}(\epsilon_\alpha+|\nabla_g u_\alpha|^2)^{\alpha-1}\big|\frac{\partial
    u_\alpha}{\partial r}\big|^2ds_0-\frac{1}{2\alpha}\dint_{\partial B_t}(\epsilon_\alpha
    +|\nabla_g u_\alpha|^2)^{\alpha-1}\dis\big|u_{\alpha,\theta}\big|^2ds_0\\[\mv]
    \s\s\s\s=\dis\frac{(\alpha-1)}{\alpha t}\dint_{B_t}(\epsilon_\alpha+|\nabla_gu_\alpha|^2
    )^{\alpha-1}|\nabla_0u_\alpha|^2dx
    +O(t).
\end{array}
\end{equation}

\subsection{Pohozaev identity}
Denote $\Delta_0=\frac{\partial^2}{\partial (x^1)^2}+\frac{\partial^2}{\partial (x^2)^2}$. By (\ref{eq1}),
we have the equation:
\begin{equation*}
\Delta_0u_\alpha+(\alpha-1)\frac{\nabla_0|\nabla_gu_\alpha|^2\nabla_0u_\alpha}
{\epsilon_\alpha+
|\nabla_gu_\alpha|^2}+A(u_\alpha)(du_\alpha,du_\alpha)=0.
\end{equation*}
As in \cite{L-W}, we multiply the both sides of the above equation
with $r\frac{\partial u_\alpha}{\partial r}$ to obtain
$$\dint_{B_t}r\frac{\partial u_\alpha}{\partial r}\Delta_0u_\alpha dx=
-(\alpha-1)\dint_{B_t}\frac{\nabla_0|\nabla_gu_\alpha|^2\nabla_0u}{\epsilon_\alpha+
|\nabla_gu_\alpha|^2}r\frac{\partial u_\alpha}{\partial r}dx.$$ It
is easy to see
$$\dint_{B_t}r\frac{\partial u_\alpha}{\partial r}\Delta_0u_\alpha dx
       =\dint_{\partial B_t}r|\frac{\partial u_\alpha}{\partial r}|^2ds_0
       -\dint_{B_t}\nabla_0(r\frac{\partial u_\alpha}{\partial r})\nabla_0u_\alpha dx.$$
Since
 $$\begin{array}{lll}
 \dint_{B_t}\nabla_0(r\frac{\partial u_\alpha}{\partial r})\nabla_0u_\alpha dx
 &=&\dint_{B_t}\nabla_0\big(x^k\frac{\partial u_\alpha}{\partial x^k}\big)\nabla_0u_\alpha dx\\[\mv]
 &=&\dint_{B_t}|\nabla_0u_\alpha|^2dx+\dint_0^t\dint_0^{2\pi}\frac{r}{2}
 \frac{\partial(|\nabla_0u_\alpha|^2)}{\partial r}rd\theta dr\\[\mv]
 &=&\dint_{B_t}|\nabla_0u_\alpha|^2dx+\frac{1}{2}\dint_{\partial B_t}|\nabla_0
 u_\alpha|^2tds_0-\dint_{B_t}|\nabla_0u_\alpha|^2dx\\[\mv]
 &=&\displaystyle{\frac{1}{2}}\dint_{\partial B_t}|\nabla_0u_\alpha|^2tds_0,
 \end{array}$$
then, we have
\begin{equation}\label{P0}
\dint_{\partial B_t}(|\frac{\partial u_\alpha}{\partial r}|^2
-\frac{1}{2}|\nabla_0
u_\alpha|^2)ds_0=-\frac{\alpha-1}{t}\dint_{B_t}
\frac{\nabla_0|\nabla_gu_\alpha|^2\nabla_0u_\alpha}{\epsilon_\alpha
+|\nabla_g u_\alpha|^2}r\frac{\partial u_\alpha}{\partial r}dx.
\end{equation}
Hence, it follows
\begin{equation}\label{P}
\dint_{\partial B_t}(|\frac{\partial u_\alpha}{\partial r}|^2
-|u_{\alpha,\theta}|^2)ds_0=-\frac{2(\alpha-1)}{t}\dint_{B_t}
\frac{\nabla_0|\nabla_gu_\alpha|^2\nabla_0u_\alpha}{\epsilon_\alpha
+|\nabla_g u_\alpha|^2}r\frac{\partial u_\alpha}{\partial r}dx.
\end{equation}
Thus, we obtain two key variational identities (\ref{e11}) and
(\ref{P}) which will be used repeatedly in our following argument.

\section{The proof of Theorem \ref{main}}

In this section, we discuss the weak energy identity on a sequence
of $\alpha$-harmonic maps. By following the idea of Ding and Tian in
\cite{D-T} we will apply (\ref{e11}) (\ref{P0}) to give the proof of Theorem
\ref{main}.

Let $B_{2\sigma}=B_{2\sigma}(0)$ be a ball in $\mathbb{R}^2$ with
the metric $g=e^{\varphi_\alpha(x)}(dx^1\otimes dx^1+dx^2\otimes dx^2)$,
where $\varphi\in C^\infty(\overline{B_{2\sigma}})$ and
$\varphi_\alpha(0)=0$, and $\varphi_\alpha$
converges smoothly. We set $u_\alpha: B_{2\sigma}\rightarrow N$ be a map
which satisfies equation (\ref{eq1}). Clearly, (\ref{e11}),
(\ref{P0}) and (\ref{P}) hold.

We assume that for any $\alpha$
$$E_{\alpha,\epsilon_\alpha}(u_\alpha,B_\sigma)<C_1, $$
and $0$ is the only blow-up point in $B_{2\sigma}$. Without loss of
generality, we assume $u_\alpha\rightarrow u_0$ in
$C^k_{loc}(B_\sigma\setminus\{0\})$, where $u_0$ is a harmonic map
from $B_\sigma$ to $N$.

We can get the first bubble in the following way. Let $x_\alpha^1\in
B_\delta$ s.t.  $|\nabla
u_{\alpha}(x_{\alpha}^1)|=\max\limits_{B_{\delta}}|\nabla
u_{\alpha}|$, and
$\lambda_{\alpha}^1=\frac{1}{\max_{B_{\delta}}|\nabla u_{\alpha}|}$.
Then, without loss of generality, we may assume in
$C^k_{loc}(\mathbb{R}^2)$
 $$u_{\alpha}(x_{\alpha}^1+\lambda_{\alpha}^1 x) \rightarrow
 v^1.$$

Now, we assume there exists  another $n_0-1$ bubbles $v^2$, $\cdots$,
$v^{n_0}$, and sequences $x_\alpha^i$, $\lambda_\alpha^i$ s.t.
$$u_\alpha(x_\alpha^i+\lambda_\alpha^ix)\rightarrow v^i$$
in $C^k_{loc}(\mathbb{R}^2\setminus A^i)$, where $A^i$ are finite
sets. Clearly, we may assume
$$\lambda_\alpha^1=\min_{i\in\{1,\cdots, n_0\}}\{\lambda_\alpha^i\}.$$
Moreover, we assume that for any $i\neq j$, one of the {\bf H1} and
{\bf H2} holds.

\subsection{The weak energy identity for the case of only one bubble}

First we prove the Theorem \ref{main} in the case of $n_0=1$, where
$n_0$ is the number of the bubbles. The general case will be
explained in the next subsection.

We denote $\lambda_\alpha=\lambda_\alpha^1$, $x_\alpha=x_\alpha^1$,
and $v=v^1$. We define
$$\Lambda_\alpha(R)=\dint_{B_{R\lambda_\alpha}(x_\alpha)}|\nabla_gu_\alpha|^{2\alpha}
dV_g,\s\Lambda=\lim_{R\rightarrow+\infty}\lim_{\alpha\rightarrow 1}\Lambda_\alpha(R).$$
and
$$\mu=\lim\limits_{\alpha\rightarrow 1}\lambda_\alpha^{2-2\alpha}.$$
By (\ref{dm}), we have $\Lambda=\mu E(v)$. Moreover, we also have
\begin{equation}\label{Lambda}
\begin{array}{ll}
&\lim\limits_{R\rightarrow+\infty}
\lim\limits_{\alpha\rightarrow 1}\dint_{B_{R\lambda_\alpha}}(\epsilon_\alpha+|\nabla_gu_\alpha|^2)^{\alpha-1}
|\nabla_gu_\alpha|^2dV_g\\[\mv]
=&
\lim\limits_{R\rightarrow+\infty}
\lim\limits_{\alpha\rightarrow 1}\dint_{B_{R}}(\epsilon_\alpha\lambda_\alpha^2+|\nabla_gv_\alpha|^2)^{\alpha-1}
\lambda_\alpha^{2-2\alpha}|\nabla_0v_\alpha|^2dx\\[\mv]
=&\mu\dint_{\mathbb{R}^2}|\nabla_0v|^2dx=\Lambda.
\end{array}
\end{equation}

Furthermore, we claim that for any $\epsilon>0$ there exist
$\delta_1$ and $R$ such that, $\forall\lambda\in
(\frac{R\lambda_\alpha}{2},4\delta_1)$, there holds
\begin{equation}\label{eqepsilon}
\int_{B_{2\lambda}\setminus B_{\lambda}(x_\alpha)}
|\nabla_gu_\alpha|^2dV_g\leq \epsilon.
\end{equation}
Suppose that the claim is false, then we may assume that there exist
$\alpha_i\rightarrow 1$ and $\lambda_i'\rightarrow 0$ satisfying
$\frac{\lambda_i'}{\lambda_{\alpha_i}}\rightarrow+\infty$ such that
\begin{equation}\label{ge}
\int_{B_{2\lambda_i'}\setminus B_{\lambda_i'}(x_{\alpha_i})}
|\nabla_gu_{\alpha_i}|^2dV_g\geq \epsilon.
\end{equation}
Denote $v_{\alpha_i}'(x)=u_{\alpha_i}(\lambda_i'x+x_{\alpha_i})$, we
may assume $v_{\alpha_i}'\rightarrow v'$ in
$C^k_{loc}(\mathbb{R}^2\setminus(\{0\}\cup \mathcal{A}),N)$, where
$\mathcal{A}$ is a finite set which does not contain $0$. If
$\mathcal{A}=\emptyset$ then it follows from (\ref{ge}) that $v'$ is
a nonconstant harmonic sphere which is different from $v^1$. This
contradicts the assumption $n_0=1$. Next, if there exists
$x_1\in\mathcal{A}$, then, by a similar argument with that we get
$v=v^1$, we can still obtain a sequence $x_i\rightarrow x_1$,
$\tilde{\lambda}_i\rightarrow 0$, s.t.
$v_i'(x_i+\tilde{\lambda}_ix)$ converges to a harmonic map $v^2$.
Hence we get
$u_{\alpha_i}(x_{\alpha_i}+\tilde{\lambda_i}(\lambda_{\alpha_i}x+x_i))$
converges to $v^2$ strongly, and then $v^2$ is the second harmonic
map. This proves that the claim (\ref{eqepsilon}) must be true.\\

\noindent Set
$$u_\alpha^*=\frac{1}{2\pi}\int_0^{2\pi}u_\alpha(x_\alpha+re^{i\theta})d\theta.$$
One is easy to check that, for any $a<b$, the following inequality
holds true
\begin{equation}\label{ms}
\begin{array}{lll}
  \dint_{B_b\setminus B_a(x_\alpha)}|\frac{\partial u_\alpha^*}{\partial r}|^2dx&=&
  \dint_a^b\int_0^{2\pi}\left|\frac{1}{2\pi}\dint_0^{2\pi}\frac{\partial u_\alpha}{\partial r}
  d\tilde{\theta}\right|^2d\theta rdr\\[\mv]
  &\leq&\dis\frac{1}{2\pi}\dint_a^b(\int_0^{2\pi}|\frac{\partial u_\alpha}{\partial
  r}|^2d\tilde{\theta}
  \dint_0^{2\pi}d\theta) rdr\\[\mv]
  &=&\dint_a^b\int_0^{2\pi}|\frac{\partial u_\alpha}{\partial r}|^2rdrd\theta
  =\dint_{B_b\setminus B_a(x_\alpha)}|\frac{\partial u_\alpha}{\partial r}|^2dx.
  \end{array}
\end{equation}

By applying (\ref{eqepsilon}) and Sacks-Uhlenbeck
$\epsilon$-regularity theorem (Theorem \ref{epsilon}), we have the following

\begin{lem}\label{upper} For any $R\lambda_\alpha<a<b<\delta_1$, we have
$$\int_{B_b\setminus B_a(x_\alpha)}|\nabla_g^2u_\alpha|r|\nabla_gu_\alpha|dV_g\leq
C\int_{B_{4b}\setminus
B_{\frac{a}{2}}(x_\alpha)}|\nabla_gu_\alpha|^2dV_g.$$
and
$$\int_{B_b\setminus B_a(x_\alpha)}|\nabla_g^2u_\alpha|\cdot|u_\alpha-u_\alpha^*|dV_g\leq
C\int_{B_{4b}\setminus
B_{\frac{a}{2}}(x_\alpha)}|\nabla_gu_\alpha|^2dV_g,$$
where $C$ does not rely on $\alpha$.
\end{lem}

\proof First, we prove the first inequality in the above lemma. We
assume that $2^Ka\in (b,2b)$ and set
$$D_i=B_{2^ia}\setminus B_{2^{i-1}a}(x_\alpha).$$

We rescale $D_i$ to $B_2\setminus B_1$, and $u_\alpha$ to
$\tilde{u}_\alpha$. By Theorem \ref{epsilon} (the $\epsilon$-regularity
theory), we have on $D_i$
$$\begin{array}{lll}
|\nabla_g u_\alpha|&\leq& \dis\frac{1}{2^{i-1}a}|\nabla_0
\tilde{u}_\alpha|_{C^0(B_2\setminus B_1)}\leq
\dis\frac{C_1}{2^{i-1}a}\|\nabla_0\tilde{u}_\alpha\|_{L^2(B_4\setminus
B_{1/2})}\\[\mv]
&& =\dis\frac{C_1}{2^{i-1}a}\|\nabla_g u_\alpha\|_{L^2(D_{i+1}\cup
D_i\cup D_{i-1})}.
\end{array}$$
Hence, it follows
\begin{equation*}
\|r\nabla_gu_\alpha\|_{C^0(D_i)}\leq 2C_1|\nabla_g
u_\alpha|\leq C_2\|\nabla_g u_\alpha\|_{L^2(D_{i+1}\cup D_i\cup
D_{i-1})}
\end{equation*}
Similarly, we have
\begin{equation*}
\|r^2\nabla_g^2u_\alpha\|_{C^0(D_i)}\leq  C_2'\|\nabla_g
u_\alpha\|_{L^2(D_{i+1}\cup D_i\cup D_{i-1})}.
\end{equation*}
Then we have
$$\begin{array}{lll}
   \dint_{D_i}|\nabla^2_gu_\alpha|r|\nabla_gu_\alpha|dV_g&\leq& C\dint_{D_{i+1}\cup
      D_i\cup D_{i-1}}|\nabla_gu_\alpha|^2dV_g\int_{D_i}\frac{dV_g}{r^2}\\[\mv]
   &\leq&C'\dint_{D_{i+1}\cup
      D_i\cup D_{i-1}}|\nabla_gu_\alpha|^2dV_g.
\end{array}$$
Therefore, we get the first inequality in this Lemma. The
proof of the second inequality goes to almost the same.
\endproof

\subsubsection{The estimate of $\int_{B_\delta\setminus B_{R\lambda_\alpha}(x_\alpha)}
| u_{\alpha,\theta}|^2dx$}

The goal of this subsection is to prove the following
\begin{lem}\label{4.1t}
For $\alpha$-harmonic map sequence $u_\alpha$ ($\alpha\rightarrow
1$), there holds true
$$\lim_{\delta\rightarrow 0}\lim_{R\rightarrow+\infty}
\lim_{\alpha\rightarrow 1}\dint_{B_\delta\setminus
B_{R\lambda_\alpha}(x_\alpha)} | u_{\alpha,\theta}|^2dx=0.$$
\end{lem}

\proof We adopt the technique of Sacks-Uhlenbeck \cite{S-U} and
\cite{Li-W} to show the lemma. Using (\ref{eqepsilon}) we have
\begin{equation}\label{osc}
   |u^*_\alpha(r)-u_\alpha(r,\theta)|\leq\epsilon_1.
\end{equation}

We compute

\begin{equation}\begin{array}{l}\label{ep1}
    \dint_{B_\delta\setminus B_{R\lambda_\alpha}(x_\alpha)}|\nabla_g u_\alpha|^2dV_g\\[\mv]
    \begin{array}{lll}
     &=&
      {\dint_{B_\delta\setminus B_{R\lambda_\alpha}(x_\alpha)}}\nabla_0 u_\alpha\nabla_0
         (u_\alpha-u^*_\alpha)dx
          +\dint_{B_\delta\setminus B_{R\lambda_\alpha}(x_\alpha)}\nabla_g u_\alpha\nabla_g
          u^*_\alpha dV_g\\[\mv]
      &=&-\dint_{B_\delta\setminus B_{R\lambda_\alpha}(x_\alpha)}\Delta_0 u_\alpha (u_\alpha-u^*_\alpha)dx
          +\dint_{B_\delta\setminus B_{R\lambda_\alpha}(x_\alpha)}\nabla_0 u_\alpha\nabla_0
          u^*_\alpha dx\\[\mv]
      & &+\dint_{\partial (B_\delta\setminus B_{R\lambda_\alpha}(x_\alpha))}\frac{\partial u_\alpha}{\partial r}
        (u_\alpha-u_\alpha^*)ds_0\\[\mv]
      &=&\dint_{B_\delta\setminus B_{R\lambda_\alpha}(x_\alpha)}A(u_\alpha)(du_\alpha,du_\alpha)
         (u_\alpha-u_\alpha^*)dV_g\\[\mv]
      &&+(\alpha-1)\dint_{B_\delta\setminus B_{R\lambda_\alpha}(x_\alpha)}\frac{\nabla_g|\nabla_g u_\alpha
         |^2\nabla_gu_\alpha}{\epsilon_\alpha+|\nabla_gu_\alpha|^2}(u_\alpha-u_\alpha^*)
          dV_g\\[\mv]
      &&+\dint_{\partial(B_\delta\setminus B_{R\lambda_\alpha}(x_\alpha))}\frac{\partial u_\alpha}{\partial r}
        (u_\alpha-u_\alpha^*)ds_0+\dint_{B_\delta\setminus B_{R\lambda_\alpha}(x_\alpha)}\frac{\partial u_\alpha}
        {\partial r}\frac{\partial u_\alpha^*}
        {\partial r}dx.
  \end{array}
\end{array}\end{equation}
On the other hand, noting (\ref{ms}) we have
\begin{equation}\label{4.2}
\begin{array}{lll}
    \dint_{B_\delta\setminus B_{R\lambda_\alpha}(x_\alpha)}\frac{\partial u_\alpha}
        {\partial r}\frac{\partial u_\alpha^*}
        {\partial r}dx&\leq&
        \sqrt{\dint_{B_\delta\setminus B_{R\lambda_\alpha}(x_\alpha)}
        |\frac{\partial u_\alpha}{\partial r}|^2dx\dint_{B_\delta\setminus
        B_{R\lambda_\alpha}(x_\alpha)}|\frac{\partial u^*_\alpha}{\partial r}|^2dx}\\[\mv]
        &\leq&\dint_{B_\delta\setminus B_{R\lambda_\alpha}(x_\alpha)}
        |\frac{\partial u_\alpha}{\partial r}|^2dx.
\end{array}
\end{equation}
Hence, by using Lemma \ref{upper}, (\ref{osc}), (\ref{4.2}) and
noting the following fact
$$|\frac{\nabla_g|\nabla_gu_\alpha|^2\nabla_gu_\alpha}{\epsilon_\alpha
+|\nabla_gu_\alpha|^2}|\leq |\nabla_g^2u_\alpha|,$$ we can infer
from (\ref{ep1})
$$\begin{array}{lll}
     \dint_{B_\delta\setminus B_{R\lambda_\alpha}(x_\alpha)}|\nabla_0 u_\alpha|^2
        dx&\leq& \dint_{B_\delta\setminus B_{R\lambda_\alpha}(x_\alpha)}|\frac{\partial u_\alpha}{\partial r}|^2
          dx+3C(\alpha-1)\dint_{B_{4\delta}}|\nabla_gu_\alpha|^2dV_g\\[\mv]
    &&+\dint_{\partial B_{\delta}(x_\alpha)}
        \frac{\partial u_\alpha}{\partial r}
        (u_\alpha-u_\alpha^*)ds_0\\
    &&-\dint_{\partial B_{R\lambda_\alpha}(x_\alpha)}
        \frac{\partial u_\alpha}{\partial r}
        (u_\alpha-u_\alpha^*)ds_0\\[\mv]
     &&+\epsilon_1'\dint_{B_\delta\setminus B_{R\lambda_\alpha}(x_\alpha)}|\nabla_0u_\alpha|^2dx,
  \end{array}$$
where $\epsilon_1'=\epsilon_1\|A\|_{L^\infty(M)}$.\\

Since $|\nabla_0u_\alpha|^2=|\frac{\partial u_\alpha}{\partial
r}|^2+|u_{\alpha,\theta}|^2$, we get
$$\begin{array}{lll}
     \dint_{B_\delta\setminus
        B_{R\lambda_\alpha}(x_\alpha)}|u_{\alpha,\theta}|^2
        dx&\leq& -\dint_{\partial B_{\delta}(x_\alpha)}
        \frac{\partial u_\alpha}{\partial r}
        (u_\alpha-u_\alpha^*)ds_0+\dint_{\partial B_{R\lambda_\alpha}(x_\alpha)}
        \frac{\partial u_\alpha}{\partial r}
        (u_\alpha-u_\alpha^*)ds_0\\[\mv]
     &&+C'((\alpha-1)+\epsilon).
  \end{array}$$
Keeping (\ref{eqepsilon}) in mind, we have

$$\lim_{\delta\rightarrow 0}\lim_{\alpha\rightarrow 1}\dint_{\partial B_{\delta}(x_\alpha)}
        \frac{\partial u_\alpha}{\partial r}
        (u_\alpha-u_\alpha^*)ds_0= 0,$$
and
$$\lim_{R\rightarrow+\infty}\lim_{\alpha\rightarrow 1}\dint_{\partial B_{R\lambda_\alpha}(x_\alpha)}
\frac{\partial u_\alpha}{\partial r}(u_\alpha-u_\alpha^*)ds_0 = 0.$$
Hence, we can see the above inequality implies the conclusion of
Lemma \ref{4.1t}.

\endproof

Immediately we infer from Lemma \ref{a1} and Lemma \ref{4.1t}
\begin{cor}\label{c1}
There holds true
$$\lim_{\delta\rightarrow 0}\lim_{R\rightarrow+\infty}
\lim_{\alpha\rightarrow 1}\dint_{B_\delta\setminus
B_{R\lambda_\alpha}(x_\alpha)}
(\epsilon_\alpha+|\nabla_gu_\alpha|^2)^{\alpha-1}|u_{\alpha,\theta}|^2dx=0.$$
\end{cor}

\subsubsection{The energy of the neck}
We set
$$F_\alpha(t)=\dint_{B_{\lambda_\alpha^t}(x_\alpha)}(\epsilon_\alpha+|\nabla_g
u_\alpha|^2)^{\alpha-1} |\nabla_0u_\alpha|^2dx,$$
$$E_{r,\alpha}(t)=\dint_{B_{\lambda_\alpha^t}\setminus
B_{\lambda_\alpha^{t_0}}(x_\alpha)}(\epsilon_\alpha+|\nabla_gu_\alpha|^2)^{\alpha-1}
|\frac{\partial u_\alpha}{\partial r}|^2dx,$$ and
$$E_{\theta,\alpha}(t)=\dint_{B_{\lambda_\alpha^t}\setminus
B_{\lambda_\alpha^{t_0}}(x_\alpha)}(\epsilon_\alpha+|\nabla_gu_\alpha|^2)^{\alpha-1}
|u_{\alpha,\theta}|^2dx.$$

By (\ref{e1}), for  $t\in [\epsilon,t_0]$, we have
$$(1-\frac{1}{2\alpha})E_{r,\alpha}'-\frac{1}{2\alpha}E_{\theta,\alpha}'=
\frac{\alpha-1}{\alpha}\log\lambda_\alpha
F_\alpha(t)+O(\lambda_\alpha^t\log \lambda_\alpha).$$
Then
$$(1-\frac{1}{2\alpha})E_{r,\alpha}(t)-\frac{1}{2\alpha}E_{\theta,\alpha}(t)=
\frac{1}{2}\dint_{t_0}^t[\frac{1}{\alpha}\log\lambda_\alpha^{2(\alpha-1)}
F_\alpha(t)+O(\lambda_\alpha^t\log\lambda_\alpha)]dt.$$ It is easy
to check that the sequences
$\{(1-\frac{1}{2\alpha})E_{r,\alpha}(t)-\frac{1}{2\alpha}E_{\theta,\alpha}(t)\}$
and $\{F_\alpha(t)\}$ are compact in $C^0([\epsilon,t_0])$ topology
for any $\epsilon>0$. Therefore, there exist two functions $F$ and
$E_r$ which belong to $C^0([\epsilon,t_0])$ such that, as
$\alpha\rightarrow 1$,
$$F_\alpha\rightarrow F,\s\s E_{r,\alpha}\rightarrow E_r\s\mbox{in}\s C^0([\epsilon,t_0]).$$
Hence, we infer from the above integration equality
$$E_r(t)=-\log\mu\dint_{t_0}^tFdt=-\log\mu\dint_{t_0}^t(E_r(t)+F(t_0))dt.$$
This implies that $E_r(t)\in C^1$ and
$$E_r'=-\log\mu(E_r+F(t_0)).$$
It follows
$$E_r(t)=\mu^{t_0-t}F(t_0)-F(t_0).$$\\

Next, we prove that
\begin{equation}\label{e4}
\lim\limits_{t_0\rightarrow
1}F(t_0)=\Lambda.
\end{equation}
Integrating (\ref{e1}) with respect to $t$ on the interval
$[R\lambda_\alpha, \lambda_\alpha^{t_0}]$ we obtain
$$\begin{array}{lll}
&&F_\alpha(t_0)-\dint_{B_{R\lambda_\alpha}(x_\alpha)}(\epsilon_\alpha+|\nabla_g
u_\alpha|^2)^{\alpha-1} |\nabla_0u_\alpha|^2dx\\[\mv]
&\leq& C\dint_{B_{\lambda_\alpha^{t_0}}\setminus
B_{R\lambda_\alpha}(x_\alpha)}(\epsilon_\alpha+|\nabla_g
u_\alpha|^2)^{\alpha-1} |u_{\alpha,\theta}|^2dx\\[\mv]
&&+C\dint_{R\lambda_\alpha}^{\lambda_\alpha^{t_0}}\frac{\alpha-1}{r}dr+
C(\lambda_\alpha^{t_0}-R\lambda_\alpha)
\end{array}$$
Noting the following holds true (from Corollary \ref{c1})
$$\lim_{t_0\rightarrow 1}\lim_{R\rightarrow+\infty}
\lim_{\alpha\rightarrow 1}\dint_{B_{\lambda_\alpha^{t_0}}\setminus
B_{R\lambda_\alpha}(x_\alpha)}
(\epsilon_\alpha+|\nabla_gu_\alpha|^2)^{\alpha-1}|u_{\alpha,\theta}|^2dx=0,$$
and
$$\lim_{t_0\rightarrow 1}\lim_{R\rightarrow+\infty}
\lim_{\alpha\rightarrow1}\dint_{R\lambda_\alpha}^{\lambda_\alpha^{t_0}}\frac{\alpha-1}{r}dr
= \lim_{t_0\rightarrow 1}\dis\frac{(1-t_0)}{2}\log\mu=0.$$ Thus,
(\ref{e4}) follows from the above inequality in view of (\ref{Lambda}).

On the other hand side, we have
$$\dint_{B_{\lambda_\alpha^t}(x_\alpha)}(\epsilon_\alpha+|\nabla_gu_\alpha|^2)^{\alpha-1}
|\nabla_gu_\alpha|^2dV_g = E_{r,\alpha}(t) + E_{\theta,\alpha}(t) +
F_\alpha(t_0).$$
Noting Corollary \ref{c1}, i.e. $\lim\limits_{\alpha\rightarrow
1}E_{\theta,\alpha}(t) =0$, we can deduce the following
$$\lim_{\alpha\rightarrow 1}\dint_{B_{\lambda_\alpha^t}(x_\alpha)}(\epsilon_\alpha+|\nabla_gu_\alpha|^2)^{\alpha-1}
|\nabla_gu_\alpha|^2dV_g = E_r(t) + F(t_0)=\mu^{t_0-t}F(t_0).$$
Thus, we have shown the following

\begin{lem}\label{33}
For any $t\in (0,1)$ and $\epsilon_\alpha>0$ with $\lim\limits_{\alpha
\rightarrow 1}\epsilon_\alpha^{\alpha-1}\geq \beta_0$, there holds
true
$$\lim_{\alpha\rightarrow 1}\dint_{B_{\lambda_\alpha^t}(x_\alpha)}(\epsilon_\alpha+|\nabla_gu_\alpha|^2)^{\alpha-1}
|\nabla_gu_\alpha|^2dV_g=\mu^{1-t}\Lambda.$$
\end{lem}

\noindent\textbf{{Proof of Theorem \ref{main}}} Here we restrict us to the
case of one bubble. By taking almost the same argument as we proved
(\ref{e4}), we obtain
\begin{equation}
\lim\limits_{\delta\rightarrow 0}\lim\limits_{t\rightarrow
0}\lim\limits_{\alpha\rightarrow1} \dint_{B_\delta\setminus
B_{\lambda_\alpha^t}(x_\alpha)} (\epsilon_\alpha+|\nabla_gu_\alpha|^2)^{\alpha-1}
|\nabla_gu_\alpha|^2dV_g =0,\end{equation} which leads to

$$\lim\limits_{\delta\rightarrow 0}\lim_{\alpha\rightarrow 1}\dint_{B_\delta(x_\alpha)}
(\epsilon_\alpha+|\nabla_gu_\alpha|^2)^{\alpha-1}
|\nabla_gu_\alpha|^2dV_g=\mu\Lambda=\mu^2E(v).$$
Obviously, this implies the required conclusion. So we have
completed the proof of Theorem \ref{main} in the case that $n_0=1$.

\subsection{The weak energy identity for the case of several bubbles}

For the general case that $n_0>1$, the proof can be completed by
induction in $n_0$, the number of bubbles.
\\
We set
$$\lambda'_\alpha=\max_{i}\{|x_\alpha^i-x_\alpha^1|+|\lambda_\alpha^i|\}.$$
Without loss of generality, we assume $x_\alpha^1\equiv 0$ and
$\lambda_\alpha'$ is attained by the $n_0$-th bubble, i.e.
$$\lambda'_\alpha=|x_\alpha^{n_0}|+|\lambda_\alpha^{n_0}|.$$

Let $v_\alpha(x)=u_\alpha(\lambda'_\alpha x)$. Then $v_\alpha$ will
converges to $v_0$ except finite points. Since $\lambda_\alpha^{n_0}$
and $\lambda_\alpha^1$ satisfies {\bf H1} or {\bf H2}, then we have
$\frac{|x_\alpha^{n_0}|}{\lambda_\alpha^1}\rightarrow +\infty$, or
$\frac{\lambda_\alpha^{n_0}}{\lambda_\alpha^1}\rightarrow+\infty$, and
therefore we have
$\frac{\lambda_\alpha'}{\lambda_\alpha^1}\rightarrow+\infty$. So, it
is easy to check that $0$ is a blowup point of the sequence
$\{v_\alpha\}$.

Similar to the proof of (\ref{eqepsilon}), we have for any $\epsilon>0$,
there are $\delta_1$ and $R$ s.t.
\begin{equation}\label{e2}
\int_{B_{2\lambda}\setminus B_{\lambda}(x_\alpha)}
|\nabla_gu_\alpha|^2dV_g\leq \epsilon,\s\forall\lambda\in (R\lambda_\alpha',\delta_1).
\end{equation}
We set $v_\alpha(x)=u_\alpha(\lambda_\alpha'x)$ and assume
$v_\alpha\rightharpoondown v_0$. Then using the arguments in the
above subsection (in this case $F(t_0)\rightarrow \lim\limits_{R
\rightarrow+\infty}\lim\limits_{\alpha\rightarrow
1}E_\alpha(v_\alpha,B_R)$ as $t_0\rightarrow 1$), we have
$$\lim_{\alpha\rightarrow 1}\dint_{B_\delta(x_\alpha)}(\epsilon_\alpha+|\nabla u_\alpha|^2)^{\alpha-1}
|\nabla u_\alpha|^2dV_g
=\lim_{\alpha\rightarrow 1}(\lambda_\alpha')^{2(2-2\alpha)}\lim_{R\rightarrow+\infty}
\lim_{\alpha\rightarrow 1}\dint_{B_R}(\lambda_\alpha^{'2}\epsilon_\alpha+|\nabla v_\alpha|^2)^{\alpha-1}
|\nabla v_\alpha|^2.$$
Moreover, (\ref{e2}) implies that all the blowup points lie in $B_R$ for some $R>0$.

The rest of the proof will be divided into two cases: i) $v_0$ is a non-trivial harmonic map. ii)
$v_0$ is trivial.\\

In case i), $v_0$ is a bubble, then we can assume $v_0$ is in fact
one of $v^i$'s for $i\in\{2,\cdots,n_0\}$. We set $v^{m_0}$ to be
equivalent to $v_0$, then $\lim\limits_{\alpha\rightarrow
1}(\lambda')^{2-2\alpha}=\mu_{m_0}$, and $E(v_0)=E(v^{m_0})$. Since
there is only $n_0-1$ bubbles of the sequence $\{v_\alpha\}$, by
induction, we have
$$\lim_{\alpha\rightarrow 1}\dint_{B_R}(\lambda_\alpha^{'2}\epsilon_\alpha
+|\nabla v_\alpha|^2)^{\alpha-1}
|\nabla v_\alpha|^2dV_g
=E(v_{0},B_R)+\sum_{i\neq m_0}(\frac{\mu_i}{\mu_{m_0}})^2E(v^i).$$

In case ii), one is easy to check that
$\frac{|x_\alpha^{n_0}|}{\lambda_\alpha^{n_0}}\rightarrow+\infty$. Then
$x_0=\lim\limits_{\alpha\rightarrow 1}
\frac{x_\alpha^{n_0}}{\lambda_\alpha'}$ which lies on $\partial B_1$ is
a blow-up point. Then there are at least two blowup points $0$ and
$x_0$. So,  at any blowup point of $v_\alpha$, there are most ${n_0}-1$
bubbles, and then we can use the induction. Hence, we will get
$$\lim_{\alpha\rightarrow 1}\dint_{B_R}(\lambda_\alpha^{'2}\epsilon_\alpha
+|\nabla v_\alpha|^2)^{\alpha-1} |\nabla v_\alpha|^2dV_g
=\sum_{i=1}^{n_0}(\frac{\mu_i}{\lim\limits_{\alpha\rightarrow 1}
(\lambda_\alpha')^{2-2\alpha}})^2E(v^i).$$ Thus, we complete the
proof of Theorem \ref{main}.

\section{Description and further analysis of the necks}
In this section, we always assume there is only one bubble on some
small ball $B_\delta$.

\subsection{The proof of Theorem \ref{main2} in the case $\nu=1$}
In this subsection, we assume $\nu=1$. Then we have $\mu=1$, and
\begin{equation}\label{ei}
\lim_{\delta\rightarrow 0}\lim_{R\rightarrow+\infty}\lim_{\alpha\rightarrow 1}
\dint_{B_\delta\setminus B_{R\lambda_\alpha}(x_\alpha)}|\nabla_g u_\alpha|^2dV_g
=0.
\end{equation}
We will use the arguments of Ding \cite{D2}.

For simplicity, we assume
$$P=\frac{\log\delta-\log R\lambda_\alpha}{\log 2}$$
is an integer. For any integer $k\in[1,P-1]$, we set
$$Q_{k}(t)=B_{2^{k+t}R\lambda_\alpha}\setminus
B_{2^{k-t}R\lambda_\alpha}(x_\alpha),$$
where $t+k\leq P$ and $k-t\geq 0$.

Using the same approximate method as in Section 3.1, we can conclude
that on $Q_{k}(t)$ the following inequality holds
\begin{equation}\label{www}
\begin{array}{lll}
    \dint_{Q_k(t)}|\nabla_0 u_\alpha|^2dx
      &\leq&\dint_{Q_k(t)}A(u_\alpha)(du_\alpha,du_\alpha)
         (u_\alpha-u_\alpha^*)dx\\[\mv]
      &&+C(\alpha-1)\dint_{Q_k(t+2)}|\nabla_0 u_\alpha|^2dx\\[\mv]
      &&+\dint_{\partial Q_k(t)}\frac{\partial u_\alpha}
         {\partial r}(u_\alpha-u_\alpha^*)ds_0 + \dint_{Q_k(t)}|\frac{\partial u_\alpha}
         {\partial r}|^2dx.
  \end{array}
\end{equation}
Next, we will apply Pohozaev identity (\ref{P0}) to controll the
last term in the above inequality, i.e.
$\int_{Q_k(t)}|\frac{\partial u_\alpha}{\partial r}|^2dx$. For the
sake of convenience, we set
$$H(r)=-\dint_{B_r(x_\alpha)}\frac{\nabla_g|\nabla_g u_\alpha|^2\nabla_g u_\alpha}{\epsilon_\alpha
+|\nabla_g u_\alpha|^2}r\frac{\partial u_\alpha}{\partial r}dV_g
=-\dint_{B_r(x_\alpha)}\frac{\nabla_0|\nabla_g u_\alpha|^2\nabla_0 u_\alpha}{\epsilon_\alpha
+|\nabla_g u_\alpha|^2}r\frac{\partial u_\alpha}{\partial r}dx.$$
Using Lemma \ref{upper}, we have
$$\begin{array}{lll}
|H(r)|&\leq&
\dint_{B_r\setminus B_{R\lambda_\alpha}(x_\alpha)}|\nabla^2_g u_\alpha|r|\frac{\partial u_\alpha}{\partial r}|dx
+|H(R\lambda_\alpha)|\\[\mv]
&\leq& C\dint_{B_{4\delta}(x_\alpha)}|\nabla_0 u_\alpha|^2dx+
|H(R\lambda_\alpha)|<C',
\end{array}$$
where we use the fact
$$\lim_{\alpha\rightarrow 1}|H(R\lambda_\alpha)|\leq
\lim_{\alpha\rightarrow
1}\dint_{B_{R\lambda_\alpha}(x_\alpha)}|\nabla_g^2u_\alpha|r|\nabla_gu_\alpha|
=\dint_{B_{R}}|\nabla^2_0v_\alpha|r|\nabla_0 v_\alpha|dx<C(R).$$
Therefore, combining these with (\ref{P0}) we obtain

$$\dint_{Q_k(t)}\big|\frac{\partial u_\alpha}{\partial r}\big|^2dx
-\frac{1}{2}\dint_{Q_k(t)}\big|\nabla_0 u\big|^2dx\leq C\dint_{2^{k-t}R\lambda_\alpha}^{
2^{k+t}R\lambda_\alpha} \frac{\alpha-1}{r}dr \leq C(\alpha-1)t.$$
It follows (\ref{www}) and the above inequality
\begin{equation}\label{ode}
(\displaystyle{\frac{1}{2}}-\epsilon_1)\dint_{Q_k(t)}|\nabla_0
   u_\alpha|^2dx\leq C(\alpha-1)(t+1)+\dint_{\partial Q_k(t)}\frac{\partial u_\alpha}
   {\partial r}(u_\alpha-u_\alpha^*)ds_0.
\end{equation}

On the other hand, we have
\begin{equation}\label{5.2}
\begin{array}{lll}
  \Big|\dint_{\partial Q_k(t)}\frac{\partial u_\alpha}{\partial r}
  (u_\alpha-u_\alpha^*)ds_0\Big|
  &\leq&\sqrt{\dint_{\partial Q_k(t)}|\frac{\partial u_\alpha}{\partial r}|^2ds_0
  \dint_{\partial Q_k(t)}| u_\alpha-u_\alpha^*|^2ds_0}\\[\mv]
  &\leq&\sqrt{\dint_{\partial Q_k(t)}|\dis\frac{\partial u_\alpha}{\partial r}|^2
  ds_0\dint_{\partial Q_k(t)}|u_{\alpha,\theta}|^2r^2ds_0}\\[\mv]
  &\leq&\dis\frac{1}{2}\left[\dint_{\partial Q_k(t)}|\frac{\partial u_\alpha}{\partial r}|^2r
  ds_0 + \dint_{\partial Q_k(t)}|u_{\alpha,\theta}|^2rds_0\right]\\[\mv]
  &=&\dis\frac{1}{2}\dint_{\partial Q_k(t)}r|\nabla u_\alpha|^2ds_0\\[\mv]
  &=&2^{t+k-1}R\lambda_\alpha\dint_{\partial B_{2^{t+k}R
  \lambda_\alpha}(x_\alpha)}|\nabla_0 u_\alpha|^2ds_0\\[\mv]
  & &-2^{k-t-1}R\lambda_\alpha\dint_{\partial B_{2^{k-t}R
  \lambda_\alpha}(x_\alpha)}|\nabla_0 u_\alpha|^2ds_0.
 \end{array}
\end{equation}
Let
$$f_k(t)=\int_{Q_{k}(t)}|\nabla u_\alpha|^2dx.$$
From (\ref{5.2}) we know
 $$\dint_{\partial Q_k(t)}\frac{\partial u_\alpha}{\partial r}
 (u_\alpha-u_\alpha^*)ds_0\leq\frac{1}{2\log2}f_k'(t).$$
Hence, by combining (\ref{ode}) and the above inequality we have
 $$(1-2\epsilon_1)f_k(t)\leq \frac{1}{\log2}f_k'(t)+
 C(\alpha-1)(t+1).$$
Multiplying the two sides of the above inequality by
$2^{-(1-2\epsilon_1)t}$ and integrating we obtain
$$f_k(1)\leq C2^{-(1-2\epsilon_1)t_1}f_k(t_1)+C(\alpha-1).$$
It is easy to check that, if we set
$$t_1=L_k=\left\{\begin{array}{ll}
k&if\s 2k-1\leq P\\
P-k& if\s 2k-1>P
\end{array}\right.$$
then, we get
$$\sqrt{E(u_\alpha,Q_k(1))}\leq C2^{-aL_k}\sqrt{E(u_\alpha,B_\delta\setminus B_{R\lambda_\alpha}(x_\alpha))}
 +C\sqrt{\alpha-1}$$
for some positive $a$ and $C$.

By the standard $L^p$ estimate, we have
 $$osc_{B_{2^{k+1}R\lambda_\alpha}\setminus
 B_{2^{k-1}R\lambda_\alpha}(x_\alpha)}u_\alpha\leq
 C2^{-aL_k}\sqrt{E(u_\alpha,B_\delta\setminus B_{R\lambda_\alpha}(x_\alpha))}
 +C\sqrt{\alpha-1}.$$
These inequalities imply
$$\begin{array}{lll}
osc_{B_\delta\setminus B_{R\lambda_\alpha}(x_\alpha)}u_\alpha&\leq&
C\sqrt{E(u_\alpha,B_\delta\setminus B_{R\lambda_\alpha}(x_\alpha))}\sum 2^{-aL_k}+C\sqrt{\alpha-1}P\\[\mv]
&\leq& C\sqrt{E(u_\alpha,B_\delta\setminus B_{R\lambda_\alpha}(x_\alpha))}+
C(R,\delta)\sqrt{\alpha-1}+C\log\lambda_\alpha^{-\sqrt{\alpha-1}}.
\end{array}$$
Letting $\alpha\rightarrow 1$, and then $R\rightarrow+\infty$, $\delta\rightarrow 0$, we get
$$osc_{B_\delta\setminus B_{R\lambda_\alpha}(x_\alpha)}u_\alpha\rightarrow 0.$$
Thus we proved Theorem \ref{main2} in the case $\nu=1$.

\subsection{The details of the neck when $\nu>1$}
The goal of this section is to show the neck converges to a geodesic
in $N$ and furthermore calculate the length of the geodesic.

For this sake, we will consider the behaviors of $u_\alpha$ on
$\partial B_{\lambda_\alpha^t}(x_\alpha)$ with $t\in [t_2,t_1]$,
where $0<t_2<t_1<1$. By the arguments in section 3.1.2, we can see
easily that
 $$\int_{B_{\lambda^t_\alpha}(x_\alpha)}|\nabla_gu_\alpha|^2dV_g\rightarrow
 \mu^{2-t}E(v^1)$$
in $C^0([t_2, t_1])$. Then, it is easy to yield
 $$\int_{B_{2\lambda_\alpha^t}\setminus B_{\frac{1}{2}\lambda_\alpha^t}
 (x_\alpha)}|\nabla_g u_\alpha|^2dV_g\rightarrow 0$$
in $C^0([t_2, t_1])$. Therefore, for any  $t \in [t_2,t_1]$, we have
\begin{equation}\label{osc20}
osc_{\partial B_{\lambda_\alpha^t}(x_\alpha)}u_\alpha\leq C\int_{B_{2\lambda_\alpha^t}
\setminus B_{\frac{1}{2}\lambda_\alpha^t}(x_\alpha)}|\nabla_g u_\alpha|^2dV_g\rightarrow 0,
\end{equation}
i.e. $u_\alpha|_{\partial B_{\lambda_\alpha^t}}(x_\alpha)$ will
subconverge to a point belonging to $N$. Especially, we have that,
as $\alpha\rightarrow 1$,
$$u_\alpha(\partial B_{\lambda_\alpha^{t_1}})\rightarrow y_1\in N
\s\mbox{and}\s u_\alpha(\partial B_{\lambda_\alpha^{t_2}})
\rightarrow y_2\in N.$$ \\

For simplicity, we will use ``$(r,\theta)$" to denote ``$x_\alpha+
r(\cos\theta,\sin\theta)$". Now we can state the main results of
this subsection as follows:

\begin{pro}\label{p1}
When $\nu>1$ and $0< t_2\leq t_\alpha\leq t_1<1$, we have, after passing to a subsequence,
\begin{equation}\label{l4.1}
\lim_{\alpha\rightarrow 1}\frac{1}{\alpha-1}\dint_{B_{R\lambda_\alpha^{t_\alpha}}
\setminus B_{\frac{1}{R}
       \lambda_\alpha^{t_\alpha}}(x_\alpha)}|u_{\alpha,\theta}|^2dx=0
\end{equation}
for any $R>0$, and
$$\frac{1}{\sqrt{\alpha-1}}\Big(u_\alpha(\lambda_\alpha^{t_\alpha}r,\theta)
-u_\alpha(\lambda_\alpha^{t_\alpha},0)\Big)\rightarrow \vec{a}\log
r$$ strongly in $C^k(S^1\times[\frac{1}{R},R],\mathbb{R}^n)$, where
$\theta$ is the angle parameter of the ball centered at $x_\alpha$,
$\vec{a}\in T_yN$ is a vector in $\mathbb{R}^n$ with
$$|\vec{a}|=\mu^{1-\lim\limits_{\alpha\rightarrow 1}t_\alpha}\sqrt{\frac{E(v)}{\pi}},$$
and $y=\lim\limits_{\alpha\rightarrow 1}u_\alpha(\lambda_\alpha^{t_\alpha},\theta)$.
\end{pro}

To prove Proposition \ref{p1}, we first prove the following

\begin{lem}\label{l411}
When $\nu>1$ and $0< t_2\leq t_\alpha\leq t_1<1$, we have
\begin{equation}\label{l4.2}
\lim_{\alpha\rightarrow 1}\frac{1}{\alpha-1}\dint_{B_{R\lambda_\alpha^{t_\alpha}}
\setminus B_{\frac{1}{R}
       \lambda_\alpha^{t_\alpha}}(x_\alpha)}|u_{\alpha,\theta}|^2dx<C
\end{equation}
where $C$ does not depend on $R$.
\end{lem}

\proof We set
$$Q(t)=B_{2^t\lambda_\alpha^{t_\alpha}}(x_\alpha)
\setminus B_{2^{-t}\lambda_\alpha^{t_\alpha}}(x_\alpha).$$ Here we
assume $2^t\leq\lambda_\alpha^{- \epsilon}$, where
$$\epsilon<\min\{t_2,1-t_1\}.$$

Applying (\ref{P0}), we get from (\ref{www}) the following
\begin{equation}\label{5.1}
\begin{array}{lll}
   (\displaystyle{\frac{1}{2}}-\epsilon_1)\dint_{Q(t)}|\nabla_0
   u_\alpha|^2dx
      &\leq& (\alpha-1)\big(\dint_{2^{-t}\lambda_\alpha^{t_\alpha}
      }^{2^{t}\lambda_\alpha^{t_\alpha}}
      \frac{1}{r}H(r)dr+C\dint_{B_{2^{t+2}\lambda_\alpha^{t_\alpha}}
      \setminus B_{2^{-t-1}\lambda_\alpha^{t_\alpha}}(x_\alpha)}|\nabla_0 u_\alpha|^2dx\big)\\[\mv]
   &&-\dint_{\partial Q(t)}\frac{\partial u_\alpha}
       {\partial r}
        (u_\alpha-u_\alpha^*)ds.
\end{array}
\end{equation}
For any $r\in
[\lambda_\alpha^{t_\alpha+\epsilon},\lambda_\alpha^{t_\alpha-\epsilon}]$,
it is easy to check that
$$|H(r)-H(\lambda_\alpha^{t_\alpha})|\leq
\dint_{B_{\lambda_\alpha^{t_\alpha-\epsilon}}\setminus
B_{\lambda_\alpha^{t_\alpha+\epsilon}}
(x_\alpha)}
|\nabla^2_gu_\alpha||r\frac{\partial u_\alpha}{\partial r}|dx.$$
Using Lemma \ref{upper}, we can get
 $$\dint_{B_{\lambda_\alpha^{t_\alpha-\epsilon}}\setminus
 B_{\lambda_\alpha^{t_\alpha+\epsilon}}(x_\alpha)}
 |\nabla^2_gu_\alpha||r\frac{\partial u_\alpha}{\partial r}|dx\leq
 C\dint_{B_{2\lambda_\alpha^{t_\alpha-\epsilon}}\setminus
 B_{\frac{1}{2}\lambda_\alpha^{t_\alpha+\epsilon}} (x_\alpha)}
 |\nabla_0 u_\alpha|^2dx.$$
By integrating (\ref{e1}) we obtain
\begin{equation}
\begin{array}{ll}\label{non}
    &(1-\dis\frac{1}{2\alpha})\dint^{2\lambda_\alpha^{t_\alpha-\epsilon}}
    _{\frac{1}{2}\lambda_\alpha^{t_\alpha+\epsilon}}ds
    \dint_{\partial B_s(x_\alpha)}(\epsilon_\alpha+|\nabla_g u_\alpha|^2)^{\alpha-1}\big|\frac{\partial
    u_\alpha}{\partial r}\big|^2ds_0\\[\mv]
    &-\dis\frac{1}{2\alpha}\dint^{2\lambda_\alpha^{t_\alpha-\epsilon}}
    _{\frac{1}{2}\lambda_\alpha^{t_\alpha+\epsilon}}ds
    \dint_{\partial B_s(x_\alpha)}(\epsilon_\alpha+|\nabla_g u_\alpha|^2)^{\alpha-1}\dis\frac{1}{s^2}
    \big|\frac{\partial u_\alpha}{\partial \theta}\big|^2ds_0\\[\mv]
    =&\dint^{2\lambda_\alpha^{t_\alpha-\epsilon}}_{\frac{1}{2}\lambda_\alpha^{t_\alpha+\epsilon}}
    \left(\dis\frac{(\alpha-1)}{\alpha s}\dint_{B_s(x_\alpha)}(\epsilon_\alpha+|\nabla_gu_\alpha|^2
    )^{\alpha-1}|\nabla_0u_\alpha|^2dx\right)ds
    +O(\lambda_\alpha^{2(t_\alpha-\epsilon)}).
\end{array}
\end{equation}
By Corollary \ref{c1} we know that the second term on the left hand
side of the above inequality vanishes as $\alpha\rightarrow 0$. On
the other hand, noting the fact $(\epsilon_\alpha
+|\nabla_gu_\alpha|^2)^{\alpha-1}$ is bounded, we have
\begin{equation*}
\begin{array}{ll}
&\dint^{2\lambda_\alpha^{t_\alpha-\epsilon}}_{\frac{1}{2}\lambda_\alpha^{t_\alpha+\epsilon}}
\left(\dis\frac{(\alpha-1)}{\alpha
s}\dint_{B_s(x_\alpha)}(\epsilon_\alpha+|\nabla_gu_\alpha|^2)^{\alpha-1}|\nabla_0u_\alpha|^2dx\right)ds\\
\leq & C
\dint^{2\lambda_\alpha^{t_\alpha-\epsilon}}_{\frac{1}{2}\lambda_\alpha^{t_\alpha+\epsilon}}
\dis\frac{(\alpha-1)}{\alpha s}ds=
\frac{C\epsilon}{\alpha}\log\lambda_\alpha^{2-2\alpha}.
\end{array}
\end{equation*}
Noting the fact $\lim_{\alpha \rightarrow
1}\epsilon_\alpha^{\alpha-1}=\beta_0>0$, we can infer from
(\ref{non}) that, as $\alpha$ is close to $1$ enough,
$$\dint_{B_{2\lambda_\alpha^{t_\alpha-\epsilon}}
 \setminus B_{\frac{1}{2}\lambda_\alpha^{t_\alpha+\epsilon}}(x_\alpha)}|\nabla_0
 u_\alpha|^2 dx \leq C\epsilon(\log\lambda_\alpha^{2-2\alpha}+1) +
 O(\lambda_\alpha^{2(t_\alpha-\epsilon)}).$$
So, when $\alpha$ is close to $1$ enough we can always choose
$\epsilon$ such that
 $$C\dint_{B_{2\lambda_\alpha^{t_\alpha-\epsilon}}
 \setminus B_{\frac{1}{2}\lambda_\alpha^{t_\alpha+\epsilon}}(x_\alpha)}|\nabla_0
 u_\alpha|^2 dx\leq C\epsilon(\log\mu + 1)<\epsilon_1.$$
Hence, we have
\begin{equation}\label{H}
H(\lambda_\alpha^{t_\alpha})-\epsilon_1\leq
H(r)\leq H(\lambda_\alpha^{t_\alpha})+\epsilon_1.
\end{equation}

Let $$f(t)=\int_{Q(t)}|\nabla_g u_\alpha|^2dV_g=\int_{Q(t)}|\nabla_0 u_\alpha|^2dx.$$ By using a similar
estimate with (\ref{5.2}) and (\ref{H}), we infer that as $\alpha$
is close to $1$ enough there holds
$$(1-2\epsilon_1)f(t)\leq \frac{1}{\log2}f'(t)+
(\alpha-1)(at+\epsilon_1),$$
where
$$a=4\log2H(\lambda_\alpha^{t_\alpha})+\epsilon_1.$$
Then, it is easy to see
$$(2^{-(1-2\epsilon_1)t}f)'\geq -(\alpha-1)(at+\epsilon_1)2^{-(1-2\epsilon_1)t}\log 2.$$
Hence, we get
 $$f(t)\leq 2^{-(1-2\epsilon_1)(\tau-t)}f(\tau)+\frac{\alpha-1}{1-2\epsilon_1}
 \big(\epsilon_1+at+\frac{a}{\log2}-a\tau 2^{-(1-2\epsilon_1)(\tau-t)}
 -a\frac{2^{-(1-2\epsilon_1)(\tau-t)}}
 {(1-2\epsilon_1)\log 2}\big).$$
Then, it follows
 $$f(k)\leq C_1(k)2^{-(1-2\epsilon_1)\tau}f(\tau)+\frac{\alpha-1}{1-2\epsilon_1}(\epsilon_1+ak+
 \frac{a}{\log 2} + aC_2(k)a\tau 2^{-(1-2\epsilon_1)\tau}).$$

Let $2^{\tau}=\lambda_\alpha^{-\epsilon}$. Then
\begin{equation}\label{ne1}
\begin{array}{lll}
    \dint_{B_{2^k\lambda_\alpha^{t_\alpha}}\setminus B_{\frac{1}{2^k}
       \lambda_\alpha^{t_\alpha}}(x_\alpha)}|\nabla_0 u_\alpha|^2dx
       &\leq&C(k)\lambda_\alpha^{{\epsilon}(1-2\epsilon_1)}
          +\dis\frac{\alpha-1}{1-2\epsilon_1}(H(\lambda_\alpha^{t_\alpha})4k\log 2+\frac{a}
          {\log 2}\\
           &&+C(k)\lambda_\alpha^{\epsilon(1-2\epsilon_1)}\log\lambda_\alpha).
  \end{array}
\end{equation}
On the other hand, by (\ref{P}) and (\ref{H}), we get
\begin{equation}\label{nor}
\dint_{B_{2^k\lambda_\alpha^{t_\alpha}}\setminus B_{\frac{1}{2^k}
       \lambda_\alpha^{t_\alpha}}(x_\alpha)}(|\frac{\partial u_\alpha}{\partial r}|^2-|u_{\alpha,\theta}|^2)dx
       \geq (\alpha-1)4k\log2(H(\lambda_\alpha^{t_\alpha})-\epsilon_1).
\end{equation}
Therefore, subtracting (\ref{ne1}) by (\ref{nor}) we obtain
\begin{equation}\label{nor1}
\begin{array}{lll}
    2\dint_{B_{2^k\lambda_\alpha^{t_\alpha}}\setminus B_{\frac{1}{2^k}
       \lambda_\alpha^{t_\alpha}}(x_\alpha)}|u_{\alpha,\theta}|^2dx
       &\leq&C(k)\lambda_\alpha^{{\epsilon}(1-2\epsilon_1)}
          +\dis\frac{(\alpha-1)}{1-2\epsilon_1}(2\epsilon_1 H(\lambda_\alpha^{t_\alpha})4k\log 2
          +\frac{a}{\log 2}\\
           &&+C(k)\lambda_\alpha^{\epsilon(1-2\epsilon_1)}\log\lambda_\alpha)
           +\epsilon_1(\alpha-1)4k\log2.
  \end{array}
\end{equation}
Since $\nu=\lim\limits_{\alpha\rightarrow 1}
\lambda_\alpha^{-\sqrt{\alpha-1}}>1$, it is easy to see that, for
any $m>0$,
\begin{equation}\label{order}
\lambda_\alpha^{{\epsilon}(1-2\epsilon_1)}=o((\alpha-1)^m).
\end{equation}
Then, noting (\ref{order}) and letting $\epsilon_1\rightarrow 0$ in
the above inequality (\ref{nor1}), we get

\begin{equation*}\label{t2}
\lim_{\alpha\rightarrow 1}\frac{1}{\alpha-1}\dint_{B_{2^k\lambda_\alpha^{t_\alpha}}
\setminus B_{\frac{1}{2^k}
       \lambda_\alpha^{t_\alpha}}(x_\alpha)}|u_{\alpha,\theta}|^2dx\leq \frac{a'}{2\log 2},
\end{equation*}
where $a'$ is a constant which does not depend on $R$. Thus, we
finish the proof of the lemma.
\endproof

Now, we are in the position to give the proof of Proposition \ref{p1}.

\proof First we show (\ref{l4.1}). Since lemma \ref{l411} says
 $$\dint_{B_{2^k\lambda_\alpha^{t_\alpha}}
 \setminus B_{\frac{1}{2^k}\lambda_\alpha^{t_\alpha}}(x_\alpha)}
 \frac{|u_{\alpha,\theta}|^2}{\alpha-1}dx=
 \dint_{2^{-k}\lambda_\alpha^{t_\alpha}}^{2^k\lambda_\alpha^{t_\alpha}}
 \frac{1}{r}\big(\dint_0^{2\pi}\frac{|\frac{\partial u_\alpha}{\partial\theta}|^2}
 {\alpha-1} d\theta\big) dr\leq\frac{a'}{2\log 2},$$
for any $\epsilon>0$, we can always find $k_0$ which is independent
of $\alpha$, s.t. there exist $$L_\alpha\in [2^{k_0}, 2^{2k_0}]$$
such that
$$\frac{1}{\alpha-1}\dint_{\partial B_{L_\alpha\lambda_\alpha^{t_\alpha}}(x_\alpha)}
| u_{\alpha,\theta}|^2rds=\frac{1}{\alpha-1}
\dint_0^{2\pi}|\frac{\partial u_\alpha}{\partial\theta} (L_\alpha
\lambda_\alpha^{t_\alpha},\theta)|^2 d\theta<\epsilon,$$
and
$$\frac{1}{\alpha-1}\dint_{\partial B_{\frac{1}{L_\alpha}\lambda_\alpha^{t_\alpha}}(x_\alpha)}
|u_{\alpha,\theta}|^2rds=\frac{1}{\alpha-1}\dint_{0}^{2\pi}
|\frac{\partial u_\alpha}{\partial\theta}|^2d\theta<\epsilon.$$ Then
\begin{equation*}
\begin{array}{lll}
  \Big|\dint_{\partial Q(\log L_\alpha/\log 2)}\frac{\partial u_\alpha}{\partial r}
        (u_\alpha-u_\alpha^*)ds\Big|
    &\leq&\sqrt{\dint_{\partial Q(\log L_\alpha/\log 2)}r|\frac{\partial u_\alpha}
      {\partial r}|^2ds\dint_{0}^{2\pi}|\frac{\partial u_\alpha}{\partial\theta}|^2d\theta}\\[\mv]
    &\leq&\sqrt{\epsilon(\alpha-1) \dint_{\partial Q(\log L_\alpha/\log 2)}
      r|\frac{\partial u_\alpha}{\partial r}|^2
        ds}.
 \end{array}
\end{equation*}
From Lemma \ref{a1} and (\ref{e1}), we get
$$\begin{array}{lll}
  \dint_{\partial Q(\log L_\alpha/\log 2)}
      r|\frac{\partial u_\alpha}{\partial r}|^2
        ds&\leq& C\dint_{\partial Q(\log L_\alpha/\log 2)}
      r|u_{\alpha,\theta}|^2
        ds+C(\alpha-1)+C\lambda_\alpha^{t_\alpha}\\[\mv]
   &\leq& (C+\epsilon)(\alpha-1).
   \end{array}$$
By (\ref{5.1}) and (\ref{H}), we get
 \begin{equation}\label{4.11}(1-2\epsilon_1)\dint_{B_{L_\alpha\lambda_\alpha^{t_\alpha}}
 \setminus B_{\frac{1}{L_\alpha}\lambda_\alpha^{t_\alpha}}(x_\alpha)}|\nabla_0u_\alpha|^2dx
 \leq \epsilon_1(\alpha-1)+2(\alpha-1)(2H(\lambda_\alpha^{t_\alpha})\log L_\alpha +\epsilon).
 \end{equation}
Noting (\ref{H}) we can infer from (\ref{P})

 $$\frac{1}{\alpha-1}\dint_{B_{L_\alpha\lambda_\alpha^{t_\alpha}}\setminus
 B_{\frac{1}{L_\alpha}\lambda_\alpha^{t_\alpha}}(x_\alpha)}
 \big(\big|\frac{\partial u_\alpha}{\partial r}\big|^2-|u_{\alpha,\theta}|^2\big)dx
 =\dint_{\frac{1}{L_\alpha}\lambda_\alpha^{t_\alpha}}^{L_\alpha\lambda_\alpha^{t_\alpha}}
 \frac{2}{r}H(r)dr,$$
which implies that
  $$\frac{1}{\alpha-1}\dint_{B_{L_\alpha\lambda_\alpha^{t_\alpha}}\setminus
  B_{\frac{1}{L_\alpha}\lambda_\alpha^{t_\alpha}}(x_\alpha)}
  \big(\big|\frac{\partial u_\alpha}{\partial r}\big|^2-|u_{\alpha,\theta}|^2\big)dx
  \geq 4\log L_\alpha(H(\lambda_\alpha^{t_\alpha})
  -\epsilon).$$
Combining (\ref{4.11}) with the above inequality we conclude the
following inequality holds true as $\alpha$ is close to $1$
sufficiently
$$\frac{1}{\alpha-1}\dint_{B_{L_\alpha\lambda_\alpha^{t_\alpha}}\setminus
B_{\frac{1}{L_\alpha}\lambda_\alpha^{t_\alpha}}(x_\alpha)}|u_{\alpha,\theta}|^2dx
\leq \epsilon_1+C\epsilon.$$ Thus, we have shown (\ref{l4.1}).\\

\noindent Next, we turn to proving the remaining assertions of
Proposition \ref{p1}.

For $\{t_\alpha\}\subset [t_2,t_1]$, we assume
 $$u_\alpha({\partial B_{\lambda_\alpha^{t_\alpha}}})\rightarrow
 y,\s\s \alpha\rightarrow 1.$$
As $N$ is regarded as an embedded submanifold in $\mathbb{R}^K$, for
simplicity, we may assume $y=0\in N$ and $T_yN=\mathbb{R}^n$, where
$\mathbb{R}^K=\mathbb{R}^n \times\mathbb{R}^{K-n}$. We also let
$\lambda_\alpha'=\lambda_\alpha^{t_\alpha}$,
$x_\alpha'=(\lambda_\alpha',0)+x_\alpha$ and
 $$u_\alpha'(x)=u_\alpha(\lambda_\alpha'x+x_\alpha),\s v_\alpha(x)=\frac{1}{\sqrt{\alpha-1}}
 [u_\alpha(\lambda_\alpha'x+x_\alpha)-u_\alpha(x_\alpha')].$$
By (\ref{ne1}) and Theorem \ref{epsilon}, we get
$$\|\nabla u_\alpha'\|_{C^0(B_{2^k}\setminus B_{2^{-k}})}
+\|\nabla^2 u_\alpha'\|_{C^0(B_{2^k}\setminus B_{2^{-k}})}<C(k)\sqrt{\alpha-1},$$
and then
$$\|\nabla v_\alpha\|_{C^0(B_{2^k}\setminus B_{2^{-k}})}
+\|\nabla^2 v_\alpha\|_{C^0(B_{2^k}\setminus B_{2^{-k}})}<C(k).$$
Noting that $v_\alpha(1,0)=0$, we get
$$\|v_\alpha\|_{C^0(B_{2^k}\setminus B_{2^{-k}})}<C'(k).$$
Obviously, we have the equation:
$$\Delta_0v_\alpha+\sqrt{\alpha-1}(A(y)+o(1))(dv_\alpha,dv_\alpha)+
(\alpha-1)O(|\nabla^2v_\alpha|)=0,$$
hence, the sequence
$$v_\alpha\longrightarrow v_0
 \s\s\mbox{in}\s C^k_{loc}(R^2\setminus\{0\})$$
where $v_0$ satisfies
 $$\Delta_0 v_0=0\,\s with\s v_0=v_0(|x|).$$

\noindent Set
$$v=(a_1,a_2,\cdots,a_n,0,\cdots,0)\log r.$$
We deduce from (\ref{e11}) that

 $$\begin{array}{lll}
 \dint_{\partial B_t}(\epsilon_\alpha + |\nabla
 u_\alpha|^2)^{\alpha-1}|\nabla_0 v_\alpha|^2ds_0 &=&
 \dis\frac{2\alpha}{2\alpha-1}\dint_{\partial B_t}(\epsilon_\alpha +
 |\nabla u_\alpha|^2)^{\alpha-1}| v_{\alpha,\theta}|^2ds_0\\
 &&+\dis\frac{2}{(2\alpha-1)t}\dint_{B_t}(\epsilon_\alpha + |\nabla
 u_\alpha|^2)^{\alpha-1}|\nabla_0 u_\alpha|^2dx.\end{array}$$
Recalling that
 $$F_\alpha(t)=\dint_{B_{\lambda_\alpha^t}}(\epsilon_\alpha + |\nabla
 u_\alpha|^2)^{\alpha-1}|\nabla_0 u_\alpha|^2dx$$
and keeping (\ref{l4.1}) in our minds, we infer from the above
identity and Lemma \ref{a1} that, as $\alpha\rightarrow 1$,
 $$\begin{array}{lll}
 \dint_{B_{2\lambda_\alpha'}\setminus B_{\lambda_\alpha'}}(\epsilon_\alpha
 +|\nabla u_\alpha|^2)^{\alpha-1}|\nabla_0 v_\alpha|^2dx
 &=&\dis\frac{2\alpha}{2\alpha-1}\dint_{\lambda_\alpha^{t_\alpha}}^{2
 \lambda_\alpha^{t_\alpha}}\frac{1}{t}F_\alpha(\log_{\lambda_\alpha}t)dt+o(1)\\[\mv]
 &=&\dis\frac{2\alpha}{2\alpha-1}\log 2 F_\alpha(t_\alpha) + o(1)\\[\mv]
 &\rightarrow& 2\log 2 F(\lim\limits_{\alpha\rightarrow1}t_\alpha).
  \end{array}$$
On the other hand, we have that, as $\alpha\rightarrow 1$, there
holds
$$\begin{array}{lll}
\dint_{B_{2\lambda_\alpha'}\setminus B_{\lambda_\alpha'}}
(\epsilon_\alpha+|\nabla u_\alpha|^2)^{\alpha-1}|\nabla_0 v_\alpha|^2
dx&=&\dint_{B_2\setminus B_1}\left(\epsilon_\alpha+|\nabla_g
v_\alpha|^2\frac{\alpha-1}{\lambda_\alpha^{'2}}\right
)^{\alpha-1}|\nabla_0v_\alpha|^2dx\\[\mv]
&\rightarrow& 2\pi\mu^{\lim\limits_{\alpha\rightarrow 1}t_\alpha}
|a|^2\log 2.
\end{array}$$
Hence, we get
$$\lim_{\alpha\rightarrow 1}v_\alpha=(a_1,\cdots,a_n,0,\cdots,0)\log r$$
with
$$\sum_{i=1}^ma_i^2=\frac{\Lambda}{\pi}\mu^{1-2\lim\limits_{\alpha\rightarrow 1}t_\alpha}.$$
As $v: S^2\longrightarrow N$ is the corresponding only bubble, then
the above identity can be written as
$$|\vec{a}|^2=\frac{E(v,S^2)}{\pi}\mu^{2-2\lim\limits_{\alpha\rightarrow
1}t_\alpha}.$$ Thus, we complete the proof of Proposition \ref{p1}.
\endproof

\begin{cor}\label{cor43}
Let $\alpha_k$ be a sequence s.t.
$$E_{\alpha_k}(u_{\alpha_k}, B_{\lambda_{\alpha_k}^t}(x_\alpha))\rightarrow\mu^{2-t}E(v)$$
in $C^0([t_2,t_1])$ with respect to $C^0$-norm. If $\nu>1$, then
$$\int_{\lambda_{\alpha_k}^t}^{2\lambda_{\alpha_k}^t}\frac{1}{\sqrt{\alpha-1}}
|\frac{\partial u_{\alpha_k}}{\partial r}|dr\rightarrow
\log2\mu^{1-t}\sqrt{\frac{E(v)}{\pi}}$$ in $C^0([t_2,t_1])$, and
$$\frac{1}{\sqrt{\alpha-1}}(r|\frac{\partial u_{\alpha_k}}{\partial r}|)(\lambda_\alpha^t,\theta)\rightarrow
\mu^{1-t}\sqrt{\frac{E(v)}{\pi}}$$ in $C^0([t_2,t_1]\times S^1)$.
\end{cor}

\proof We need only to prove the first claim, since the proof of the
second claim is similar. If the first claim was not true, then we
assumed that there was a subsequence $\alpha_{k_i}$, $t_i\rightarrow
t_0$ s.t.
$$\Big|\int_{\lambda_{\alpha_{k_i}}^{t_i}}^{2\lambda_{\alpha_{k_i}}^{t_i}}
\frac{1}{\sqrt{\alpha-1}}|\frac{\partial u_{\alpha_{k_i}}}{\partial
r}|dr- \log2 \mu^{1-t_i}\sqrt{\frac{E(v)}{\pi}}\Big|\geq
\epsilon>0.$$

On the other hand, from the above arguments on Proposition \ref{p1} we
know that, after passing to a subsequence, there holds
$$\frac{u_{\alpha_{k_i}}(\lambda_{\alpha_{k_i}}x)-
u_{\alpha_{k_i}}(\lambda_{\alpha_{k_i}},0)}{\sqrt{\alpha-1}}\rightarrow
\vec{a}\log r,$$ with
$|\vec{a}|=\Big|\mu^{1-t_0}\sqrt{\frac{E(v)}{\pi}}\Big|$. Hence we
derive the following
$$\lim_{i\rightarrow+\infty}\int_{\lambda_{\alpha_{k_i}}^{t_i}}^{2\lambda_{\alpha_{k_i}}^{t_i}}
\frac{1}{\sqrt{\alpha-1}}|\frac{\partial u_{\alpha_{k_i}}}{\partial
r}|dr=|\vec{a}|\int_1^2\frac{1}{r}dr=\log2 \mu^{1-t_0}
\sqrt{\frac{E(v)}{\pi}}.$$ This is a contradiction.
\endproof

\subsection{The proof Theorem \ref{main2} in the case $\nu>1$}

First, we need to show the necks for the $\alpha$-harmonic map
sequence converge to some geodesics in $N$ which join the bubbles.
For this goal, we denote the curve
$$\frac{1}{2\pi} \int_0^{2\pi}
u_\alpha(r,\theta)d\theta:[\lambda_\alpha^{t_1},
\lambda_\alpha^{t_2}] \longrightarrow \mathbb{R}^n$$ by
$\Gamma_\alpha$. For simplicity, we set
$$\omega_\alpha(r)=\frac{1}{2\pi}\int_{0}^{2\pi}u_\alpha(r,\theta)d\theta.$$
First, we claim that if $\Gamma_\alpha\rightarrow\Gamma$, then $\Gamma$
must lies on $N$. This is a direct corollary of (\ref{osc20}).
Next, we will prove a subsequence of $\Gamma_\alpha$
will converges locally to a geodesic of $N$ and then give the formula of
length of $\Gamma$.

By computation we have
$$\begin{array}{lll}
  \ddot{\omega}_\alpha&=&\displaystyle\frac{1}{2\pi}\dint_{0}^{2\pi}\ddot{u}_\alpha(r,\theta)d\theta\\[\mv]
    &=&\displaystyle\frac{1}{2\pi}\dint_{0}^{2\pi}(\ddot{u}_\alpha(r,\theta)
    +\displaystyle\frac{u_{\alpha,\theta\theta}}{r^2})d\theta\\[\mv]
    &=&\displaystyle\frac{1}{2\pi}\dint_{0}^{2\pi}\Delta_0 u_\alpha d\theta
    -\displaystyle\frac{1}{2\pi}\int_{0}^{2\pi}\frac{\dot{u}_\alpha}
      {r}d\theta\\[\mv]
    &=&-\displaystyle\frac{1}{2\pi}\dint_0^{2\pi}A(u_\alpha)(du_\alpha,du_\alpha)-\frac{\alpha-1}{2\pi}
     \dint_0^{2\pi}\frac{\nabla_0|\nabla_g u_\alpha|^2\nabla_0 u_\alpha}{\epsilon_\alpha^2+|\nabla_g u_\alpha|^2}
      d\theta-\frac{\dot\omega_\alpha}{r}
\end{array}$$
where we have used the fact
$$\int_0^{2\pi}u_{\theta\theta}(r,\theta)d\theta=0.$$

Let
$$G_\alpha = -\ddot{\omega}_\alpha-\frac{\dot{\omega}_\alpha}{r}.$$
Denote the induced metric of $\Gamma_\alpha$ in $\mathbb{R}^K$ by
$h_\alpha$, and let $A_{\Gamma_\alpha}$ be the second fundamental
form of $\Gamma_\alpha$ in $\mathbb{R}^K$.

Given
$\lambda_\alpha^{t_\alpha}\in[\lambda_\alpha^{t_1},\lambda_\alpha^{t_2}]$.
As before, we always have
$$\frac{u_\alpha(\lambda_\alpha^{t_\alpha}r,
 \theta)-u_\alpha(\lambda_\alpha^{t_\alpha},0)}{\sqrt{\alpha-1}} \rightarrow \vec{a}\log r,$$
where $\vec{a}\in T_yN$ and $y=\lim\limits_{\alpha\rightarrow
1}u_\alpha(\lambda_\alpha^{t_\alpha},\theta)$. Therefore, we have
\begin{equation}\label{AC}
\dot{\omega}_\alpha(\lambda_\alpha^{t_\alpha})
=\frac{\sqrt{\alpha-1}}{\lambda_\alpha^{t_\alpha}}(\vec{a}+o(1)),\s
h_\alpha(\frac{d}{dr},\frac{d}{dr})=|\dot\omega_\alpha|^2=\frac{\alpha-1}{\lambda_\alpha^{2t_\alpha}}
(|\vec{a}|^2+o(1)).
\end{equation}
where $o(1)\rightarrow 0$ as $\alpha\rightarrow 1$. Moreover, we have
$$\begin{array}{lll}
G_\alpha(\lambda_\alpha^{t_\alpha})&=&\dis\frac{1}{2\pi}\dint_0^{2\pi}A(u_\alpha)(du_\alpha,du_\alpha)d\theta
+\frac{\alpha-1}{2\pi}
     \dint_0^{2\pi}\frac{\nabla_0|\nabla_g u_\alpha|^2\nabla_0 u_\alpha}{\epsilon_\alpha^2+|\nabla_g u_\alpha|^2}
      d\theta\\[\mv]
     &=&\dis\frac{\alpha-1}{\lambda_\alpha^{2t_\alpha}}(\frac{1}{2\pi}\dint_0^{2\pi}A(y)
     (\vec{a},\vec{a})d\theta+o(1))+(\alpha-1)\dint_0^{2\pi}O(|\nabla^2_gu_\alpha|)d\theta\\[\mv]
     &=&\dis\frac{\alpha-1}{\lambda_\alpha^{2t_\alpha}}(A(y)(\vec{a},\vec{a})+o(1)+O(\sqrt{\alpha-1}))\\[\mv]
     &=&\dis\frac{\alpha-1}{\lambda_\alpha^{2t_\alpha}}(A(y)(\vec{a},\vec{a})+o(1)).
\end{array}$$
In the above identity we have made use of the fact $\nu>1$ which
implies that for any $m>0$
$$\lambda_\alpha^{2t_\alpha}=o((\alpha-1)^m).$$ Noting that
$\left\langle A(y)(\vec{a},\vec{a}), \vec{a}\right\rangle=0$, we get
\begin{equation}\label{AC2}
\begin{array}{lll}
-A_{\Gamma_\alpha}(d\omega_\alpha,d\omega_\alpha)&=&\ddot{\omega}_\alpha-
\dis\frac{\lan\ddot{\omega}_\alpha,\dot{\omega}_\alpha\ran}
{|\dot{\omega}|^2}\dot{\omega}_\alpha
=-G_\alpha+\dis\frac{\lan G_\alpha,\dot{\omega}_\alpha\ran\dot{\omega}_\alpha}{|\dot{\omega}_\alpha|^2}\\[\mv]
&=&-\dis\frac{\alpha-1}{\lambda_\alpha^{2t_\alpha}}(A(y)(\vec{a},\vec{a})+o(1)).
\end{array}
\end{equation}
Hence, we get
$$\|A_{\Gamma_\alpha}\|_{h_\alpha}^2(\lambda_\alpha^{t_\alpha})<C.$$
Similar to the proof of Corollary \ref{cor43}, we have, after passing to a subsequence,
$$\|A_{\Gamma_\alpha}\|_{h_\alpha}^2(\lambda_\alpha^{t})<C$$
for any $t\in [t_2,t_1]$.

Now, we fix $y\in N$, and let $s$ to be the arc length parameter of
$\omega_\alpha(t)$ with $s(\lambda_\alpha^{t_\alpha})=0$. We assume
$\omega_\alpha(\lambda_\alpha^{t_\alpha})\rightarrow y$ as
$\alpha\rightarrow 1$. It is well-known that
$\|A_{\Gamma_\alpha}\|_{h_\alpha}^2(\lambda_\alpha^{t_\alpha})$ does
not depend on the choice of parameter, and
$$\frac{d^2\omega_\alpha}{ds^2}= -A_{\Gamma_\alpha}(\omega_\alpha)
(\frac{d\omega_\alpha}{ds},\frac{d\omega_\alpha}{ds}),$$ then
$\omega_\alpha(s)$ will converges locally to $\omega(s)$ in $C^1$,
where $s$ is still the arc length parameter. This implies that
$\Gamma_\alpha|_{[\lambda_\alpha^{t_1},\lambda_\alpha^{t_2}]}$
converges locally to a curve $\Gamma$ locally. We claim that
$$A_{\Gamma_\alpha}(\omega_\alpha)(\frac{d\omega_\alpha}{ds},\frac{d\omega_\alpha}{ds})\rightarrow
A(\omega)(\frac{d\omega}{ds},\frac{d\omega}{ds})$$ strongly in
$C^0([0,s_1],\mathbb{R}^n)$ for sufficiently small $s_1$. If this
was not true, then for any small $s_1$ we could find a subsequence
of $\{u_\alpha\}$, still denoted by $\{u_\alpha\}$, such that

$$s_\alpha'=\int_{\lambda_\alpha^{t_\alpha}}^{\lambda_\alpha^{t_\alpha'}}|\dot\omega_\alpha|
dr\rightarrow s'\in(0,s_1)$$
s.t.
$$\Big|A_{\Gamma_\alpha}(\omega_\alpha)(\frac{d\omega_\alpha}{ds},\frac{d\omega_\alpha}{ds})-
A(\omega)(\frac{d\omega}{ds},\frac{d\omega}{ds})\Big|_{s=s_\alpha'}>\epsilon.$$
To apply Proposition \ref{p1},
we must ensure that $t_\alpha'\in[\frac{t_2}{2},t_1]$.
For simplicity, we may assume
$\lambda_\alpha^{\frac{t_2}{2}}=2^P\lambda_\alpha^{t_\alpha}$ where
$P$ is an integer. Then, applying Corollary \ref{cor43} we have
 $$\dint_{2^i\lambda_\alpha^{t_\alpha}}^{2^{i+1}\lambda_\alpha^{t_\alpha}}|\dot\omega_\alpha|dr
 =\sqrt{\alpha-1}\mu^{1-(t_\alpha-i\log_{\lambda_\alpha}2)}\left(\log 2\sqrt{\frac{E(v)}{\pi}}+o_\alpha(1)\right).$$
Therefore, as $\alpha$ is close to $1$ enough,
 $$\begin{array}{lll}\dint_{\lambda_\alpha^{t_\alpha}}^{\lambda_\alpha^{\frac{t_2}{2}}}|\dot\omega_\alpha|dr
 =\sum_{i=0}^{P-1}\dint_{2^i\lambda_\alpha^{t_\alpha}}^{2^{i+1}\lambda_\alpha^{t_\alpha}}|
 \dot\omega_\alpha|dr
 &\geq& P\sqrt{\alpha-1}\left(\sqrt{\dis\frac{E(v)}{\pi}}\log 2+o_\alpha(1)\right)\\
 &\geq & C(t_\alpha-\dis\frac{t_2}{2})\log\lambda_\alpha^{-\sqrt{\alpha-1}}\\
 &\geq & C\dis\frac{t_2}{2}\log\nu>0.\end{array}$$
Therefore, we may always choose $s_1$ to be very small, for example
$s_1<C\frac{t_2}{2}\log\nu$, such that $t_\alpha'\in [\frac{t_2}{2},
t_1]$. Then, as before there holds
 $$\frac{u_\alpha(\lambda_\alpha^{t_\alpha'}r,\theta)-
 u_\alpha(\lambda_\alpha^{t_\alpha'},0)}{\sqrt{\alpha-1}}\rightarrow
 \vec{a^\prime}\log r.$$
Obviously,
$$\frac{\dot\omega_\alpha(\lambda_\alpha^{t_\alpha'})}{|\dot\omega_\alpha(\lambda_\alpha^{t_\alpha'})|}
=\frac{d\omega_\alpha}{ds}(s_\alpha') \rightarrow
\frac{d\omega}{ds}(s').$$
Applying (\ref{AC}) and (\ref{AC2}), we
get that, after passing a subsequence the following holds
$$A_{\Gamma_\alpha}(\omega_\alpha)(\frac{d\omega_\alpha}{ds},\frac{d\omega_\alpha}{ds})|_{s=s_\alpha'}=
\frac{1}{|\dot{\omega}_\alpha(\lambda_\alpha^{t_\alpha'})|^2}A_{\Gamma_\alpha}(\omega_\alpha)
(\dot{\omega}_\alpha,\dot{\omega}_\alpha)|_{r=\lambda_\alpha^{t_\alpha'}}\rightarrow
A(\omega)(\frac{d\omega}{ds},\frac{d\omega}{ds})|_{s=s'}$$ which
contradicts the choice of $s_\alpha'$. So, we infer that
$$\frac{d\omega}{ds}(s)-\frac{d\omega}{ds}(0)= -\int_0^sA(\omega)(\frac{d\omega}{ds},\frac{d\omega}{ds})ds$$
holds near $s=0$. This shows $\omega$ is smooth near $0$ and
satisfies
 $$\frac{d^2\omega}{ds^2}= -A(\omega)(\frac{d\omega}{ds},\frac{d\omega}{ds}).$$
Therefore, we obtain finally
 $$\nabla_{\frac{d\omega}{ds}}^N\frac{d\omega}{ds}=\frac{d^2\omega}{ds^2} +
 A(\omega)(\frac{d\omega}{ds},\frac{d\omega}{ds})=0,$$
which means that $\Gamma$ is a geodesic.\\

Next, we calculate the length of the geodesic $\Gamma$.
 For simplicity, we assume $\lambda_\alpha^{t_2}=2^P\lambda_\alpha^{t_1}$
for some integer $P$. Then we have
$$P=\frac{t_2-t_1}{\log 2}\log\lambda_\alpha.$$\\

When $\nu=+\infty$,
by Corollary \ref{cor43}, we have
$$L(\Gamma_\alpha|_{B_{2^{k+1}\lambda_\alpha^{t_1}\setminus 2^k\lambda_\alpha^{t_1}}(x_\alpha)})
\geq\sqrt{\alpha-1}(\sqrt{\frac{E(v)}{\pi}}\log 2+o(1)).$$
Then
$$L(\Gamma_\alpha)\geq CP\sqrt{\alpha-1}\geq C\log\lambda_\alpha^{-\sqrt{\alpha-1}}\rightarrow+\infty.$$
This implies
$$L(\Gamma)=+\infty.$$\\

Now, we assume $\nu<+\infty$. By Corollary \ref{cor43},
$$L(\Gamma_\alpha|_{B_{2^{k+1}\lambda_\alpha^{t_1}\setminus 2^k\lambda_\alpha^{t_1}}(x_\alpha)})
=\sqrt{\alpha-1}(\sqrt{\frac{E(v)}{\pi}}\log 2+o_\alpha(1)),$$
where $o_\alpha(1)\rightarrow 0$ as $\alpha\rightarrow 1$ uniformly.
Hence
 $$L(\Gamma)=\lim_{\alpha\rightarrow 1}\sqrt{\alpha-1}P\sqrt{\frac{E(v)}{\pi}}\log 2
 =(t_1-t_2)\sqrt{\frac{E(v)}{\pi}}\log\nu.$$
Now, it is easy to see that to complete the proof of Theorem
\ref{main2} we only need to prove the following:
\begin{equation}\label{osc1}
 osc_{B_{\lambda_\alpha^t}\setminus B_{R\lambda_\alpha}(x_\alpha)}u_\alpha\rightarrow 0,\s \mbox{as}\s
 \alpha\rightarrow 1, \s\mbox{then}\s R\rightarrow+\infty\s and \s t\rightarrow 1
 \end{equation}
and
 \begin{equation}\label{osc2}
 osc_{B_{\delta}\setminus B_{\lambda_\alpha^t}(x_\alpha)}u_\alpha\rightarrow 0,\s \mbox{as}\s
 \alpha\rightarrow 1, \s\mbox{then}\s\delta\rightarrow 0\s and \s t\rightarrow 0.
 \end{equation}
Since $\nu<+\infty$ implies $\mu=1$, from Theorem \ref{main} we know
 $$\lim_{t\rightarrow 1}\lim_{R\rightarrow+\infty}\lim_{\alpha\rightarrow 1}
 \int_{B_{\lambda_\alpha^t}\setminus B_{R\lambda_\alpha}(x_\alpha)}|\nabla u_\alpha|^2=0.$$
Therefore, we can use the same method as in Subsection 4.1 (we
replace $\delta$ with $\lambda_\alpha^t$) to deduce
$$osc_{B_{\lambda_\alpha^t}\setminus B_{R\lambda_\alpha}}u_\alpha\leq
C\sqrt{E(u_\alpha,B_{\lambda_\alpha^t}\setminus
B_{R\lambda_\alpha}(x_\alpha))}+C(1-t)\log\nu+C\sqrt{\alpha-1},$$ then
(\ref{osc1}) follows. Similarly, we can prove (\ref{osc2}). Hence,
we derive the length formula of the geodesic $\Gamma$
 $$L=\sqrt{\frac{E(v)}{\pi}}\log\nu.$$
Thus, we finish the proof of Theorem \ref{main2}.
\endproof

Next, we want to give the proof of Corollary \ref{intro.c1}.
However, to prove the corollary we only need to prove the following
proposition:

\begin{pro}\label{ho} If $\nu<+\infty$, then when $(\alpha-1)$ is sufficiently small,
all the $u_\alpha$ are in the same homotopy class.
\end{pro}

\proof When $\alpha_i-1$ and $\alpha_j-1$ are sufficiently small, we
have $\|u_i-u_j\|_{C^0}\leq i(N)$ where $i(N)$ is the injective
radius of $N$. Hence, by using exponential map we know that
$u_{\alpha_i}$ and $u_{\alpha_j}$ are homotopic in $M\setminus
B_\delta$, $B_{\lambda_\alpha^{t_2}} \setminus
B_{\lambda_\alpha^{t_1}}$ and $B_{R\lambda_\alpha}$ respectively.

Let $p=u_0(0)$ and $q=v(+\infty)$. By (\ref{osc1}), we know that
$u_i(B_\delta(x_\alpha)\setminus
B_{\lambda_\alpha^{t_2}}(x_\alpha))$ is contained in a simply
connected ball centered at $p$ when $\alpha$ is close to 1 enough,
$\delta$ and $t_2$ are small enough. Similarly, by (\ref{osc2}) we
also have $u_j(B_{\lambda_\alpha^{t_1}}\setminus
B_{R\lambda_\alpha}(x_\alpha))$ is contained in a small simply
connected ball in $N$ with center $q$ when $\alpha-1$, $\delta$ and
$1-t_1$ are sufficiently small. Hence $u_i$ and $u_j$ are homotopic
in $B_{\delta}\setminus B_{\lambda_\alpha^{t_2}}$ and
$B_{\lambda_\alpha^{t_1}}\setminus B_{R\lambda_\alpha}$
respectively. So $u_i$ and $u_j$ are homotopic.
\endproof

\section{Some comments and an example}
In this paper we only consider the case $u_\alpha$ is an
$\alpha$-harmonic maps when the conformal structure of $M$ is fixed.
Naturally, one will ask the following problems (i) what could we say
in the case $u_\alpha$ is an $\alpha$-harmonic maps and the
conformal structure of $M$ varies with $\alpha$, (ii) whether the
methods in this paper can be extended to a class of variational
problem which is more general than $\alpha$-energy or not. In a
forthcoming paper we will further develop some tools to discuss some
issues which relate to the above problems.

On the other hand, one want to know whether one can give an example
to show there is a neck joining the bubbles in the limit of an
$\alpha$-harmonic map sequence is of infinite length or not.
However, if we can construct a manifold N and find a minimizing
$\alpha$-harmonic map sequence which satisfies the condition of
Corollary \ref{intro.c1}, then the corollary tells us that indeed
there exists a necks in the limit which if of infinite length. By
modifying the example of Duzaar and Kuwert (see page 304 of
\cite{D-K}) we can construct such example as following .\\

{\noindent\bf Example.} Let $\mathbb{Z}^3$ act on $\mathbb{R}^3$ by
$\tau_{\kappa}(x,y,z)= (x+4k_1,y+4k_2, z+4k_3)$, where
$\kappa=(k_1,k_2,k_3)\in\mathbb{Z}^3$. Consider
$$\tilde{X}=\mathbb{R}^3\setminus\cup_{\kappa}\tau_\kappa(B_1(0))$$
and $X$ is the quotient of $\tilde{X}$. Then $X$ is a compact
manifold with boundary. Topologically, $X$ is $T^3$ minus a small
ball.

Let $\Phi$ be a conformal map from $\mathbb{R}^2$ to
$\partial{B_1(0)}$, s.t. $\Phi(x)=(1,0,0)$ when $|x|>2$ and
$\Phi(x)=(-1,0,0)$ when $|x|<1$, and $deg(\Phi)=1$ if we consider
$\Phi$ be a map from $S^2$ to $S^2$. Moreover, we let
$\gamma_k:[0,1]\rightarrow\tilde{X}$ be a curve connect $(4k-1,0,0)$
and $(1,0,0)$. We define
$$v_k=\left\{\begin{array}{ll}
\Phi(x), &|x|\geq\delta\\
\gamma_k\left(\dis\frac{\log r-\log R\epsilon}{\log\delta-\log
R\epsilon}\right), &R\epsilon<|x|<\delta\\
\tau_{(k,0,0)}(\Phi(\frac{x}{\epsilon})).&\frac{|x|}{\epsilon}\leq R
\end{array}\right.$$
We denote $\pi$ to be the projection from $\tilde{X}$ to $X$, then
$\pi(v_k)\in\pi_2(X)$. We have
$$\begin{array}{lll}
  \dint_{B_\delta\setminus B_{R\epsilon}}|\nabla v_k|^2&=&
     2\pi\dint_{R\epsilon}^\delta|\frac{\partial\gamma_k}{\partial r}|^2rdr\\
     &<&c\dis\frac{\|\dot{\gamma}\|_{L^\infty}^2}
     {(-\log R\epsilon+\log\delta)^2}
\dint_{R\epsilon}^\delta\frac{dr}{r}
   =c\dis\frac{\|\dot{\gamma}\|_{L^\infty}^2}{\log\delta-\log{R\epsilon}},
\end{array}$$
$$\dint_{R^2\setminus B_\delta}|\nabla v_k|^2 \leq E(\Phi),\s and \s
\dint_{B_{R\epsilon}}|\nabla v_k|^2 \leq E(\Phi).$$ So, we can find
suitable $R$ and $\epsilon$, s.t.
$$E(\pi(u_k))=E(u_k)\leq 2E(\Phi)+1.$$

We claim that $[\pi(v_k)]$ are different homotopy classes. Assuming
this is not true, we can find a continuous map
$$H(x,t):S^2\times[0,1]\rightarrow X$$
s.t.
$$H(x,0)=\pi(v_i)\s and \s H(x,1)=\pi(v_j).$$
Since $S^2\times[0,1]$ is simply connected,
 we are able to lift $H$ to $\tilde{H}$
which is a map from $S^2\times[0,1]\rightarrow\tilde{X}$ with
$\tilde{H}(x,0)= v_i$. We assume that
$\tilde{H}(x,1)=\tau_\kappa(v_j)$. Hence $[v_i]=[\tau_\kappa(v_j)]$.
Therefore
$$[\partial B_1(0)+\partial \tau_{(i,0,0)}(B_1(0))]
=[\partial \tau_\kappa(B_1(0))+\partial
\tau_{(j,0,0)}\tau_\kappa(B_1(0))]\s \mbox{in}\
\s\pi_2(\tilde{X}),$$ where $\pi_2(\tilde{X})$ is the second
homotopy group of $\tilde{X}$. However, it is easy to check that
$\pi_1(\tilde{X})=\{1\}$, then by Hurewicz Theorem, the above
identity is not true.

Now, we proceed to construct $N$. Let $f$ be a homeomorphism from
$X$ to $Y=X$. We consider the quotient space of $X\cup Y$, obtained
by gluing every point $x\in\partial X$ with $f(x)\in \partial Y$
together. In this way, we get a closed compact manifold $N$ and a
projection $\phi:N\rightarrow X$. One is easy to check that
$\pi(v_k)$ can be also considered as a map from $S^2$ to $N$ with
$E(\pi(v_k))<C$. We claim that $[\pi(v_k)]$ are some different
homotopic classes with each other in $\pi_2(N)$. Assuming it is not
true. Then, we can find a continuous map $H(x,t):
S^2\times[0,1]\rightarrow N$ such that $H(x,0)=\pi(v_i)$ and
$H(x,1)=\pi(v_j)$. Hence, $\phi(H(x,t))$ is
just a homotopic map of $\pi(v_i)$ and $\pi(v_j)$ in $X$. A contradiction.\\

\noindent Finally, we would like to ask the following problems:\\

\noindent{\bf Problem 1.} Suppose all $\alpha$-harmonic maps
$u_\alpha$ belong to the same homotopic class and satisfy the energy
identity as $\alpha\rightarrow 1$. Do the necks consist of some
geodesics of finite length?\\

\noindent{\bf Problem 2.} Could we find a sequence
$\alpha_k\rightarrow 1$, and $\alpha_k$ harmonic maps
$u_{\alpha_k}$, s.t. 1)the Morse index tends to infinite; 2)
$\sup_kE_{\alpha_k}(u_{\alpha_k})<\infty$; 3) for any $i\neq j$,
$u_{\alpha_i}$ and $u_{\alpha_j}$ are not  homotopic to each other.

{}

Yuxiang Li

Department of Mathematical Sciences, Tsinghua University,

Beijing 100084, P.R.China.

Email: yxli@math.tsinghua.edu.cn; liyuxiang@tsinghua.edu.cn\\

Youde Wang

{Academy of Mathematics and Systematic Sciences, Chinese Academy of
Sciences,

Beijing 100080,  P.R. China.

Email: wyd@math.ac.cn}

\begin{thebibliography}{2}

\bibitem[C]{C} Chang, K.C.: Heat flow and boundary value problem for harmonic maps.
 Ann. Inst. H. Poincare Anal. Non Lineaire  {\bf 6}  (1989),  no. 5, 363--395.

\bibitem[C-T]{C-T} Chen, J. and Tian, G.: Compactification of moduli space
of harmonic mappings.  Comment. Math. Helv.  {\bf 74}  (1999), 201-237.

\bibitem[C-M]{C-M} Colding, T. H. and Minicozzi II, W. P.: Width and finite extinction
time of Ricci flow. Preprint.


\bibitem[D]{D2} Ding, W.: Lectures on the heat flow of harmonic maps. Manuscript.


\bibitem[D-T]{D-T} Ding, W. and Tian, G.: Energy identity for a class of approximate
harmonic maps from surfaces, Comm. Anal. Geom. {\bf 3} (1995),
543-554.

\bibitem[D-K]{D-K} Duzaar, F and Kuwert, E: Minimization of conformally invariant
energies in homotopy classes.
Calc. Var. Partial. Differ. Equ., {\bf 6} (1998), 285-313.

\bibitem[E-S]{E-S} Eells, J. and Sampson, J. H.: Harmonic mappings of Riemannian manifolds.
Amer. J. Math.  {\bf 86}  (1964), 109--160.

\bibitem[J]{J} Jost, J.: Two-dimensional geometric variantional problems. John Wiley and Sons, Chichester,
1991.

\bibitem[L]{L} Lamm, T.: Energy identy for approximations of harmonic maps from surfaces. Preprint.

\bibitem[L-W]{Li-W} Li, Y. and Wang, Y.: Bubbling location for a sequence of approximate
$f$-harmonic maps from surfaces. Preprint.

\bibitem[Lin-W]{L-W} Lin, F. and Wang, C.: Energy identity of harmonic
map flows from surfaces at finite singular time.  Calc. Var. Partial
Differential Equations  {\bf 6}  (1998),   369-380.

\bibitem[M]{M} Moore, J.D.: Energy growth in minimal surface bubbles. ``http://www.math.ucsb
.edu/~moore/growthrev.pdf''.

\bibitem[P]{P} Parker, T. H.:
Bubble tree convergence for harmonic maps.
J. D. G. {\bf 44} (1996), 595-633.

\bibitem[Q]{Q} Qing, J.:  On singularities of the heat flow for harmonic maps from surfaces into
spheres. Comm. Anal. Geom., {\bf 3} (1995), 297-315.

\bibitem[Q-T]{Q-T} Qing, J. and Tian, G.: Bubbling of the heat flows for harmonic maps from surfaces.
Comm. Pure. Apple. Math., {\bf 50} (1997), 295-310.

\bibitem[S-U]{S-U} Sacks, J. and Uhlenbeck, K.: The existence of minimal immersions of
2-spheres, {\it Ann. of Math.}, {\bf 113} (1981), 1-24.

\bibitem[St]{St} Struwe, M.: On the evolution of harmonic mappings of Riemannian surfaces,
Comment. Math. Helv. {\bf 60} (1985), 558-581.


\bibitem[T]{T} Topping, P.: Winding behaviour of finite-time singularities of the harmonic map heat flow.
Math. Z., {\bf 247} (2004), 279-302.\\
\end{thebibliography}
\end{document}